\documentclass[12pt]{amsart}
\usepackage[T2A]{fontenc}
\usepackage[english]{babel}
\usepackage{amsaddr}
\usepackage{amsmath}
\usepackage{amssymb}
\usepackage{amsfonts}
\usepackage{epsfig}
\usepackage{srcltx}
\usepackage{subfigure}
\usepackage{cite}
\usepackage{float}
\usepackage{mathtools}
\usepackage[a4paper, mag=1000, includefoot, left=2cm, right=2cm, top=2cm, bottom=2cm, headsep=1cm, footskip=1cm]{geometry}
  \usepackage[unicode,colorlinks]{hyperref}
  \hypersetup{colorlinks=true, citecolor=blue, linkcolor=blue}

\newtheorem{Th}{Theorem}
\newtheorem{Lem}{Lemma}

\begin{document}

\thispagestyle{empty}

\title[ ]{Asymptotic regimes in oscillatory systems with damped non-resonant perturbations}

\author[O.A. Sultanov]{Oskar A. Sultanov}

\address{
%Chebyshev Laboratory, St. Petersburg State University, 14th Line V.O., 29, Saint Petersburg 199178 Russia;\\
Institute of Mathematics, Ufa Federal Research Centre, Russian Academy of Sciences, Chernyshevsky street, 112, Ufa 450008 Russia.}
\email{oasultanov@gmail.com}

%\thanks{\it \today}

\maketitle

{\small
\begin{quote}
\noindent{\bf Abstract.} 
An autonomous system of ordinary differential equations describing nonlinear oscillations on the plane is considered. The influence of time-dependent perturbations decaying at infinity in time is investigated. It is assumed that the perturbations satisfy the non-resonance condition and do not vanish at the equilibrium of the limiting system. Possible long-term asymptotic regimes for perturbed solutions are described. In particular, we show that the perturbed system can behave like the corresponding limiting system or new asymptotically stable regimes may appear. The proposed analysis is based on the combination of the averaging technique and the construction of Lyapunov functions.

 \medskip

\noindent{\bf Keywords: }{asymptotically autonomous system, nonlinear oscillations, damped perturbations, asymptotics, stability, averaging}

\medskip
\noindent{\bf Mathematics Subject Classification: }{34D10, 34D20, 34D05, 37J40, 70K65}
%34D10 Perturbations of ordinary differential equations
%34D05 Asymptotic properties of solutions to ordinary differential equations
%34D20 Stability of solutions to ordinary differential equations
%37J40 Perturbations of finite-dimensional Hamiltonian systems, normal forms, small divisors, KAM theory, Arnol’d diffusion
%70K65 Averaging of perturbations for nonlinear problems in mechanics

\end{quote}
}
{\small

\section*{Introduction}
The paper is devoted to studying the effect of perturbations decaying in time on the class of autonomous systems. Such perturbed asymptotically autonomous systems arise in a wide range of problems, for example, in the study of epidemiological models~\cite{CCT95}, Painlev\'{e} equations~\cite{BG08}, phase-locking phenomena~\cite{OS18}, celestial motions~\cite{DS22}, and in other problems associated with nonlinear non-autonomous systems~\cite{KF13,Pan21,Dong22}. Note that global properties of solutions of asymptotically autonomous systems have been studied in many papers~\cite{LRS02,KS05,MR08}. In particular, it follows from~\cite{LM56} that, under some conditions, decaying perturbations can preserve the dynamics described by the limiting systems. See also~\cite{LDP74}, where the conditions are described under which damped perturbations do not disturb the global dynamics of autonomous oscillatory systems. In the general case, the global properties of perturbed and unperturbed trajectories can differ significantly~\cite{HT94,OS21IJBC}. It depends both on the class of decaying perturbations and on the class of limiting systems.

In this paper, we study the effect of decaying perturbations on nonlinear oscillatory systems. It is assumed that the perturbations oscillate with an asymptotically constant frequency, satisfy the non-resonance condition, and do not preserve the equilibrium of the limiting system. Note that similar perturbations, but vanishing at the equilibrium, have been considered in several papers. For example, linear systems were considered in \cite{HL75,MP85,PN07,BN10}, where the asymptotic behaviour of solutions at infinity was investigated. Bifurcations of the equilibrium and the asymptotic behaviour of perturbed solutions for the corresponding nonlinear systems were discussed in \cite{OS21DCDS,OS21JMS} for both the resonant and non-resonant cases. Decaying oscillatory perturbations with increasing frequency in time were studied in~\cite{DF78,BD79,OS22JMS,OS23DCDSB}. However, to the best of the author's knowledge, damped perturbations with an asymptotically constant frequency that do not preserve equilibrium have not previously been considered in detail. In this case, difficulties arise both at the stage of choosing some analogues of fixed points and at the analysis of their stability. In this paper, we consider such perturbations in the non-resonant case and discuss long-term asymptotic regimes of solutions and their stability.

The paper is organized as follows. In Section~\ref{PS}, the statement of the problem is given and the class of perturbations is described. The main results are presented in Section~\ref{MR}. The justification is provided in the subsequent sections. In particular, in Section~\ref{Sec3}, a near-identity transformation that simplifies the system in the first asymptotic terms at infinity in time is constructed. Then, under some natural assumptions on the structure and parameters of the simplified equations, possible asymptotic regimes for solutions of the perturbed system and their stability are discussed. The justification of these results for different cases is contained in sections~\ref{SQ1},~\ref{SQ2} and~\ref{SQ3}. In Section~\ref{SEx} the proposed theory is applied to some examples of asymptotically autonomous systems. The paper concludes with a brief discussion of the results obtained.

\section{Problem statement}\label{PS}
Consider a non-autonomous system of ordinary differential equations
\begin{gather}\label{FulSys}
\frac{dr}{dt}=f(r,\varphi,S(t),t), \quad 
\frac{d\varphi}{dt}=\omega(r)+r^{-1} g(r,\varphi,S(t),t),
\end{gather}
where the functions $\omega(r)$, $f(r,\varphi,S,t)$ and $g(r,\varphi,S,t)$ are infinitely differentiable, defined for all $|r|\leq r_0={\hbox{\rm const}}$, $(\varphi,S)\in\mathbb R^2$, $t>0$, and $2\pi$-periodic with respect to $\varphi$ and $S$. It is assumed that $S(t)\sim s_0 t $ as $t\to\infty$ with $s_0={\hbox{\rm const}}>0$, and for each fixed $r$ and $\varphi$
\begin{gather*}
f(r,\varphi,S(t),t)\to 0, \quad g(r,\varphi,S(t),t)\to 0, \quad t\to\infty.
\end{gather*}
Hence, system \eqref{FulSys} is asymptotically autonomous, and the limiting system is given by
\begin{gather}\label{LimSys}
\frac{dr}{dt}=0, \quad \frac{d\varphi}{dt}=\omega(r).
\end{gather} 
 We assume that $\omega(r)>0$ for all $|r|\leq r_0$. In this case, system \eqref{LimSys} describes nonlinear oscillations on the plane $(x,y)=(r\cos\varphi,-r\sin\varphi)$ with a period $T(r)\equiv 2\pi/\omega(r)$ and an equilibrium at the origin $(0,0)$. The solutions $r(t)$ and $\varphi(t)$ of system \eqref{FulSys} corresponds to the amplitude and the phase of oscillations. The functions $f(r,\varphi,S(t),t)$ and $r^{-1}g(r,\varphi,S(t),t)$ play the role of perturbations of system \eqref{LimSys}.

Note that system \eqref{FulSys} in the Cartesian coordinates $(x,y)$ takes the form 
\begin{gather*}%\label{SysXY}
\frac{dx}{dt}=\omega\left(r\right) y + F(x,y,t), \quad 
\frac{dy}{dt}=-\omega\left(r\right) x + G(x,y,t),
\end{gather*}
where $r=\sqrt{x^2+y^2}$, $F(x,y,t)\equiv (x f +y g)/r$ and $G(x,y,t)\equiv (y f -x g)/r$. We see that the perturbations $F(x,y,t)$ and $G(x,y,t)$  do not generally preserve the equilibrium $(0,0)$. 
 
The goal of the paper is to study the influence of perturbations $f(r,\varphi,S(t),t)$ and $r^{-1} g(r,\varphi,S(t),t)$ on the dynamics in the vicinity of the solution $r=0$ of the limiting system and to describe possible long-term asymptotic regimes.

Let us specify the class of perturbations. We assume that 
\begin{gather}\label{FG}\begin{split}
&f(r,\varphi,S,t)\sim \sum_{k=1}^\infty t^{-\frac{k}{q}} f_k(r,\varphi,S), \\ 
&g(r,\varphi,S,t)\sim \sum_{k=1}^\infty t^{-\frac{k}{q}} g_k(r,\varphi,S), \quad 
t\to\infty,
\end{split}
\end{gather}
for all $|r|\leq r_0$ and $(\varphi,S)\in\mathbb R^2$, where $q\in\mathbb Z_+$ and the coefficients $f_k(r,\varphi,S)$, $g_k(r,\varphi,S)$ are $2\pi$-periodic with respect to $\varphi$ and $S$ such that their Fourier series 
\begin{gather*}
f_k(r,\varphi,S)=\sum_{(k_1,k_2)\in\mathbb Z^2} f_{k,k_1,k_2}(r) e^{i (k_1 \varphi+k_2 S)}, \quad
g_k(r,\varphi,S)=\sum_{(k_1,k_2)\in\mathbb Z^2} g_{k,k_1,k_2}(r) e^{i (k_1 \varphi+k_2 S)}
\end{gather*}
contain a finite number of harmonics $e^{ik_2 S}$. Moreover, it assumed that the following condition holds:
\begin{gather} \label{pFG}
\begin{split}
\exists\, p\in\mathbb Z_+:  \quad  &|f_k(r,\varphi,S)|+ |g_k(r,\varphi,S)|=\mathcal O(r), \quad k<p, \quad  f_{p,0,0}(0)\neq 0.
\end{split}
\end{gather}
Note that this assumption ensures that the first $p-1$ terms in the series \eqref{FG} preserve the equilibrium of the limiting system.
The phase of perturbations is considered in the form
\begin{gather}\label{Sform}
S(t)=s_0 t + \sum_{k=1}^{q-1} s_k t^{1-\frac{k}{q}}+s_q \log t,
\end{gather}
where $s_k={\hbox{\rm const}}$ and the parameter $s_0$ satisfies a non-resonance condition
\begin{gather}\label{nres}
k_1\omega(0)+k_2 s_0\neq 0 \quad \forall\, k_1,k_2\in\mathbb Z, \quad |k_1|+|k_2|\neq 0.
\end{gather}

Consider the example
\begin{gather}\label{Ex0}
\frac{dx}{dt}=y, \quad 
\frac{dy}{dt}=-x+ t^{-1} \lambda y + t^{-\frac{p}{2}} \frac{\gamma(S(t))y}{\sqrt{x^2+y^2}},
\end{gather}
where $\gamma(S)=\gamma_0+\gamma_1\sin S$, $S(t)=\sqrt 2 t$, $\lambda, \gamma_0, \gamma_1\in\mathbb R$, and $p\in\mathbb Z_+$. It can easily be checked that system \eqref{Ex0} corresponds in the polar coordinates $x=r\cos\varphi$, $y=-r\sin\varphi$ to \eqref{FulSys} with 
\begin{align*}
&f(r,\varphi,S,t)\equiv t^{-1} \lambda r\sin^2\varphi+t^{-\frac p2}\gamma(S) \sin^2\varphi, \\ 
&g(r,\varphi,S,t)\equiv t^{-1}\lambda r \sin\varphi\cos\varphi+t^{-\frac p2} \gamma(S)\sin\varphi\cos\varphi,
\end{align*}
$s_0=\sqrt 2$, $\omega(r)\equiv 1$, $f_{p,0,0}(0)=\gamma_0/2$. If $\lambda=\gamma_0=\gamma_1=0$, system \eqref{Ex0} has $2\pi$-periodic general solution $x(t)\equiv \rho \cos(t+\phi)$, $y(t)\equiv-\rho \sin(t+\phi)$, where $\rho$ and $\phi$ are arbitrary constants. In this case, $r(t)\equiv \rho$ and $\varphi(t)\equiv t+\phi$. In the absence of oscillating part of the perturbation ($\lambda\neq 0$, $\gamma(S)\equiv 0$), the asymptotics of the solutions can be easily constructed by using the WKB approximations~\cite{WW66}:
\begin{gather*}
r(t)=\rho t^{\frac{\lambda}{2}}\left(1+\mathcal O(t^{-1})\right), \quad
\varphi(t)=t+\phi+\mathcal O(t^{-1}), \quad t\to\infty.
\end{gather*} 
Note that the global behaviour of the system in this case is determined by the sign of the parameter $\lambda$. Numerical analysis of the system with $\lambda\neq 0$, $\gamma_0\neq 0$ and $\gamma_1\neq 0$ shows that the global dynamics depends both on the value of the parameters $\lambda$, $\gamma_0$, $\gamma_1$ and on the degree of decay $p/2$ of the oscillatory part of the perturbation. In particular, if $p>2$, we see that, as in the previous case, the stability of the system depends on the sign of $\lambda$ (see~Fig.~\ref{Fig1}, a). However, if $p\leq 2$, oscillatory term of damped perturbations can break the stability of the system (see~Fig.~\ref{Fig1}, b, c).  

\begin{figure}
\centering
\subfigure[$p=3$]{
 \includegraphics[width=0.3\linewidth]{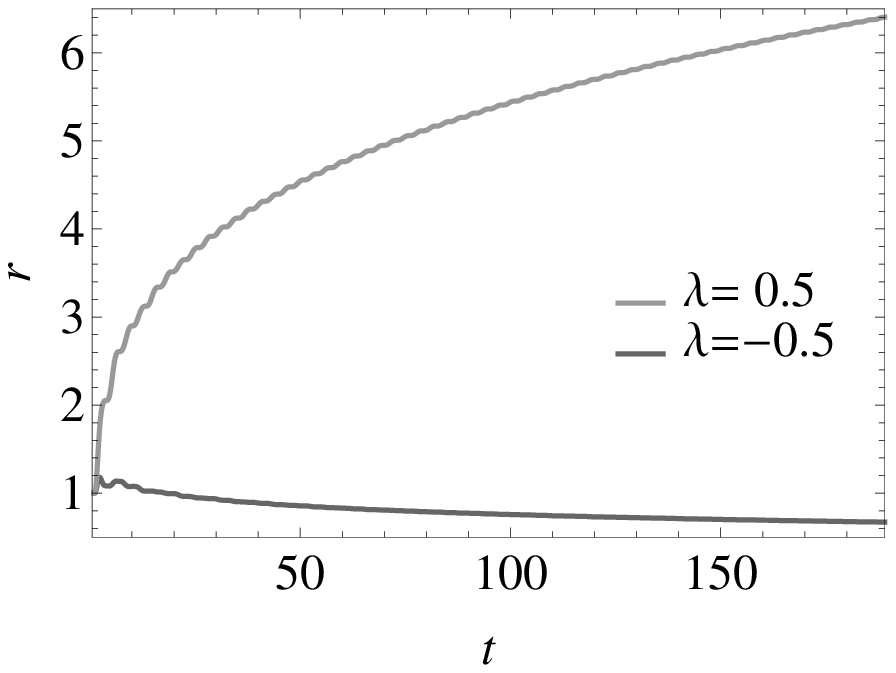}
}
\hspace{1ex}
 \subfigure[$p=2$]{
 	\includegraphics[width=0.3\linewidth]{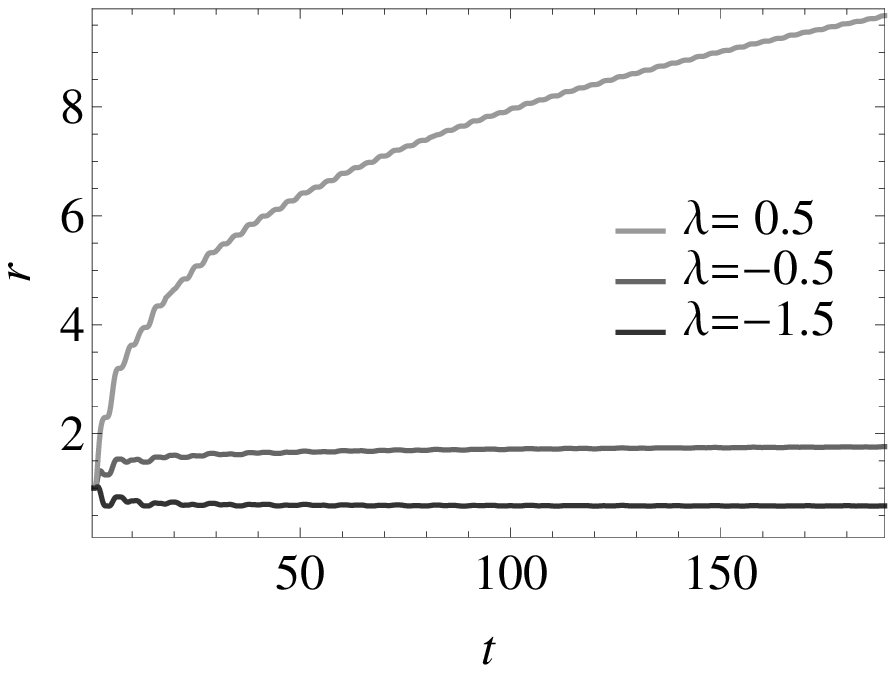}
}
\hspace{1ex}
\subfigure[$p=1$]{
 	\includegraphics[width=0.3\linewidth]{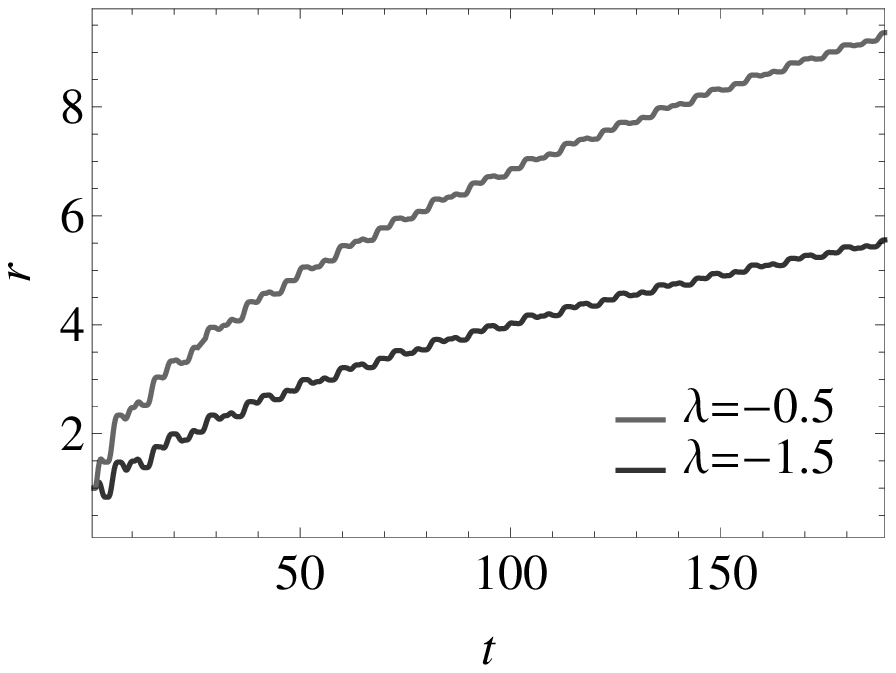}
}
\caption{\small The evolution of $r(t)\equiv \sqrt{x^2(t)+y^2(t)}$ for solutions to system \eqref{Ex0} with $\gamma_0=\gamma_1=1$ and different values of the parameters $p$ and $\lambda$.} \label{Fig1}
\end{figure}

Thus, the paper investigates the role of the parameters and the structure of the considered class of damped perturbations on the qualitative and asymptotic behaviour of solutions.

%and there exist $r _0>0$ and $r_0>0$ such that for all $r \in (0,r _0]$ the level lines $\{(x,y)\in\mathbb R^2: H(x,y)=r ^2/2\}$ lying in the ball $\mathcal B_{r_0}=\{(x,y): r\leq r_0\}$ define a family of closed curves on the phase space $(x,y)$ parametrized by the parameter $r $. To each closed curve there correspond a periodic solution $\xi_0(t,r )$, $\eta_0(t,r )$ of system \eqref{LimSys} with a period $T(r )\equiv 2\pi/\omega(r )$, where $\omega(r )\neq 0$ for all $r \in [0,r _0]$. For definiteness, assume that $\xi_0(0,r )>0$ and $\eta_0(0,r )=0$ for all $r \in (0,r _0]$. It can easily be checked that $\omega(r ) = 1 + \mathcal O(r ^2)$ as $r \to 0$. Note that the value $r  = 0$ corresponds to the isolated fixed point $(0, 0)$ of center type. Moreover, it is assumed that $\mathcal B_{r_0}$ does not contain any fixed points of the limiting system, other than the origin.

\section{Main results}\label{MR}

Let us present the main results of the paper. The first of them is devoted to averaging of the system in the leading asymptotic terms and constructing the corresponding near-identity transformation.

\begin{Th}\label{Th1}
Let system \eqref{FulSys} satisfy \eqref{FG}, \eqref{pFG}, \eqref{Sform} and \eqref{nres}.
Then for all $N\in [p,2p)$ there exist $r_\ast\in (0,r_0]$, $t_\ast>0$ and the transformation
\begin{gather}\label{rhosubs0}
\rho(t)=r (t)+\tilde v_N(r (t),\varphi(t),t), \quad  \tilde v_N(r,\varphi,t)\equiv t^{-\frac 1q}\tilde v_{N,1}(r,\varphi,t)+t^{-\frac pq} \tilde v_{N,p}(r,\varphi,t)
\end{gather}
such that for all $|r|\leq r_\ast$, $\varphi\in\mathbb R$ and $t\geq t_\ast$ system \eqref{FulSys} can be transformed into
\begin{gather}\label{FulSys3}
\frac{d\rho}{dt}= \Lambda_N(\rho,t)+\tilde \Lambda_N(\rho,\varphi,t), \quad 
\frac{d\varphi}{dt}=\omega(\rho)+\tilde \Omega_N(\rho,\varphi,t)
\end{gather}
where the functions $\Lambda_N(\rho,t)$, $\tilde \Lambda_N(\rho,\varphi,t)$, $\tilde \Omega_N(\rho,\varphi,t)\equiv \tilde \Omega_{N,1}(\rho,\varphi,t)+\rho^{-1}\tilde \Omega_{N,p}(\rho,\varphi,t)$ are defined for all $|\rho|\leq \rho_\ast$, $\varphi\in\mathbb R$, $t\geq t_\ast$, $2\pi$-periodic with respect to $\varphi$, and satisfy
\begin{gather}\nonumber
\Lambda_N(\rho,t)\equiv \sum_{k=1}^N t^{-\frac{k}{q}}\Lambda_k(\rho), \quad 
\Lambda_k(\rho)=\begin{cases} \mathcal O(|\rho|), & k<p,\\ 
\mu_{p}+\mathcal O(|\rho|), & k=p,\\ 
\mathcal O(1), & k> p, \end{cases} \quad \rho\to 0,\\
 \label{LNONest}
 \tilde \Lambda_N(\rho,\varphi,t)=\mathcal O(t^{-\frac{N+1}{q}}), \quad 
 \tilde \Omega_{N,1}(\rho,\varphi,t)=\mathcal O(t^{-\frac{1}{q}}), \quad 
\tilde \Omega_{N,p}(\rho,\varphi,t)=\mathcal O(t^{-\frac{p}{q}}), \quad t\to\infty
\end{gather}
uniformly for all $\varphi\in\mathbb R$ and $|\rho|\leq \rho_\ast$ with $\mu_{p}=f_{p,0,0}(0)\neq 0$. Moreover, $\tilde v_{N,1}(r,\varphi,t)=\mathcal O(|r|)$ and $ \tilde v_{N,p}(r,\varphi,t)=\mathcal O(1)$ as $r\to 0$ uniformly for all $\varphi\in\mathbb R$ and $t\geq t_\ast$.
\end{Th}

The proof is contained in Section~\ref{Sec3}. 

Thus, the dynamics of system \eqref{FulSys} is determined by the behaviour of solutions of the transformed equations. Consider some additional assumptions on the right-hand sides of system \eqref{FulSys3}.

First note that the transformation described in Theorem~\ref{Th1} can set some terms in $\Lambda_N(\rho,t)$ to zero. Let
$n\leq p$ be the smallest number such that
\begin{gather}\label{as1}
	 \Lambda_k(\rho)\equiv 0, \quad k\leq n-1, \quad \Lambda_n(\rho)\not\equiv 0.
\end{gather}

Along with \eqref{FulSys3}, we consider the truncated system
\begin{align}\label{trSys}
&\frac{d\varrho}{dt}=\Lambda_N(\varrho,t), \\
\label{trSys2}
&\frac{d\phi}{dt}=\omega(\varrho) 
\end{align}
as $t\geq t_\ast$. This system is obtained from \eqref{FulSys3} by dropping residual terms $\tilde \Lambda_N(\rho,\varphi,t)$ and $\tilde \Omega_N(\rho,\varphi,t)$. 

It can be shown that long-term regimes for solutions to system \eqref{FulSys} are associated with the solutions of the truncated system and depend on the values of $(n,p)$.

Consider the following three cases (see Fig.~\ref{Part}):
\begin{align*}
&\mathcal Q_1=\{(n,p)\in\mathbb Z^2_+:\quad n< p, \quad n\leq q\},\\
&\mathcal Q_2=\{(n,p)\in\mathbb Z^2_+:\quad n= p, \quad n\leq q\},\\
&\mathcal Q_3=\{(n,p)\in\mathbb Z^2_+:\quad n\leq p, \quad n> q\}.
\end{align*}

\begin{figure}
\centering
 \includegraphics[width=0.3\linewidth]{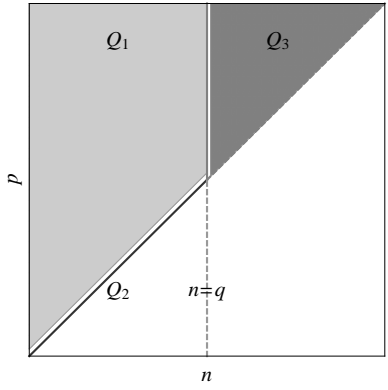}
\caption{\small Partition of the parameter plane $(n,p)$.} \label{Part}
\end{figure}

\subsection{The case $(n,p)\in\mathcal Q_1$}

Assume that $\Lambda_n(\rho)$ has a nonzero linear part in the vicinity of the equilibrium: 
\begin{gather}\label{aslin}
	\Lambda_n(\rho)=\rho (\lambda_n+\mathcal O(|\rho|)), \quad \rho \to 0
\end{gather}
with $\lambda_n={\hbox{\rm const}}\neq 0$. Define $\nu_0=(p-n)/q>0$. Then, we have

\begin{Lem}\label{Lem1}
Let $(n,p)\in\mathcal Q_1$ and assumptions \eqref{as1}, \eqref{aslin} hold with $\lambda_n+\delta_{n,q}\nu_0<0$. Then system \eqref{trSys}, \eqref{trSys2} has a solution $\varrho_1(t)$, $\phi_1(t)$ such that
\begin{gather}\label{sol11}
\varrho_1(t)\sim   
\sum_{k=0}^{\infty} t^{-\nu_0-\frac{k}{q}} \rho_k(\log t), \quad \phi_1(t)  \sim  \omega(0) t, \quad t\to\infty,
\end{gather}
 where $\rho_k(\tau)$ are some polynomials in $\tau$ and $\rho_0(\tau)\equiv -{\mu_p}/(\lambda_n+\delta_{n,q}\nu_0)$. 
\end{Lem}

Moreover, the solution $\varrho_1(t)$ is asymptotically stable. 
\begin{Lem}\label{Cor1}
Let $(n,p)\in\mathcal Q_1$ and assumptions \eqref{as1}, \eqref{aslin} hold with $\lambda_n+\delta_{n,q}\nu_0<0$. Then the solution $\varrho_1(t)$ of equation \eqref{trSys} is exponentially {\rm (}polynomially{\rm)} stable if $n<q$  $(n=q)$. 
\end{Lem}

It can be shown that the trajectories of the full system \eqref{FulSys} behave like the solution $\varrho_1(t)$, $\phi_1(t)$ of the truncated system.
\begin{Th}\label{Th2}
Let system \eqref{FulSys} satisfy \eqref{FG}, \eqref{pFG}, \eqref{Sform}, \eqref{nres}, and $(n,p)\in \mathcal Q_1$ be integers such that assumptions \eqref{as1}, \eqref{aslin} hold with $\lambda_{n}+\delta_{n,q}\nu_0<0$. Then for all $\varepsilon>0$ there exist $\delta_0>0$ and $t_0\geq t_\ast$ such that the solution $r(t)$, $\varphi(t)$ of system \eqref{FulSys} with initial data $|r(t_0)-\varrho_1(t_0)|\leq \delta_0$ and $\varphi(t_0)\in \mathbb R$ satisfies 
\begin{gather}\label{rineq1}
\sup_{t\geq t_0}\left\{t^{\nu_0}|r(t)-\varrho_1(t)|\right\}<\varepsilon,\\
\label{varphiineq1} 
\varphi(t)\sim \phi_1(t), \quad  t\to\infty.
\end{gather}
\end{Th}

On the other hand, the asymptotic regime described in Lemma~\ref{Lem1} turns out to be unstable if $\lambda_n+\delta_{n,q}\nu_0>0$. 
Let $\varrho_{1,M}(t)$ be a $M$th partial sum of the series  \eqref{sol11}. Then, we have 
\begin{Th}\label{Lem1u}
Let system \eqref{FulSys} satisfy \eqref{FG}, \eqref{pFG}, \eqref{Sform}, \eqref{nres}, and $(n,p)\in\mathcal Q_1$ be integers such that assumptions \eqref{as1}, \eqref{aslin} hold with $\lambda_n+\delta_{n,q}\nu_0>0$. Then for all $M\in\mathbb Z_+$ there exists $\varepsilon>0$ such that for all $\delta_0\in (0,\varepsilon)$ the solution  $r(t)$, $\varphi(t)$ of system \eqref{FulSys} with initial data $|r(t_s)-\varrho_{1,M}(t_s)|\leq \delta_0$, $\varphi(t_s)\in\mathbb R$ at some $t_s\geq t_\ast$  satisfies the inequality $t_e^{\nu_0}|r(t_e)-\varrho_{1,M}(t_e)|\geq \varepsilon$ at some $t_e>t_s$.
\end{Th}

Note that if $n=q$ and $\lambda_n<0$, it can be shown that solutions of system \eqref{FulSys} with small enough initial amplitude will remain close to zero. 

\begin{Th}\label{ThDop}
Let system \eqref{FulSys} satisfy \eqref{FG}, \eqref{pFG}, \eqref{Sform}, \eqref{nres}, and $(n,p)\in\mathcal Q_1$ be integers such that assumptions \eqref{as1}, \eqref{aslin} hold with $n=q$ and $\lambda_n<0$. Then for all $\varepsilon>0$ there exist $\delta_0>0$ and $t_0\geq t_\ast$ such that the solution  $r(t)$, $\varphi(t)$ of system \eqref{FulSys} with initial data $|r(t_0)|\leq \delta_0$ and $\varphi(t_0)\in\mathbb R$ satisfies the inequality $|r(t)|< \varepsilon$ for all $t\geq t_0$.
\end{Th}

Now, consider the case, when $\Lambda_n(\rho)$ is strongly nonlinear as $\rho\to 0$: 
\begin{gather}\label{asnonlin}
\exists\, d, m \in\mathbb Z_+: \quad m\geq 2, \quad 
	\Lambda_k(\rho)=
	\begin{cases}
	\rho^m (\lambda_{n,m}+\mathcal O(|\rho|)), & k=n,\\
	\mathcal O(|\rho|^{m}), & k<n+d,\\
	\rho(\lambda_{n+d}+\mathcal O(|\rho|)), & k=n+d,
	\end{cases}
\end{gather}
with $\lambda_{n,m},\lambda_{n+d}={\hbox{\rm const}}\neq 0$. 
Define the parameters
\begin{gather*}
m_\ast=\frac{p-n}{p-(n+d)}, 
\quad 
\vartheta_m=
\begin{cases}
\displaystyle \frac{p-n}{qm}, & m\leq m_\ast,\\
\displaystyle  \frac{d}{q(m-1)}, & m>m_\ast,
\end{cases}
\end{gather*}
and the function
\begin{gather*}
\mathcal P(z)\equiv 
	\begin{cases}
	\displaystyle 	\lambda_{n,m} z^m+\mu_p, & m<m_\ast,\\
	\displaystyle 	\lambda_{n,m} z^m+(\lambda_{n+d}+\delta_{n+d,q}\vartheta_m)z+\mu_p, & m=m_\ast,\\
	\displaystyle 	\lambda_{n,m} z^m+(\lambda_{n+d}+\delta_{n+d,q}\vartheta_m)z, & m>m_\ast.
	\end{cases}
\end{gather*}
Assume that
\begin{gather}\label{Pmas}
\exists\, z_m\in\mathbb R:\quad \mathcal P(z_m)=0, \quad \mathcal P'(z_m)\neq 0.
\end{gather}
Then, we have the following:

\begin{Lem}\label{LemM}
Let $(n+d,p)\in\mathcal Q_1$ and assumptions \eqref{as1}, \eqref{asnonlin}, \eqref{Pmas} hold with $\mathcal P'(z_m)<0$. Then system \eqref{trSys}, \eqref{trSys2} has a solution $\varrho_m(t)$, $\phi_m(t)$ such that
\begin{gather}\label{solm}
\varrho_m(t)\sim   
\sum_{k=0}^{\infty} t^{-\vartheta_m-\frac{k}{qm}} \rho_k(\log t), \quad \phi_m(t)  \sim  \omega(0) t, \quad t\to\infty,
\end{gather}
 where $\rho_k(\tau)$ are some polynomials in $\tau$ and $\rho_0(\tau)\equiv z_m$. 
\end{Lem}

\begin{Th}\label{ThM}
Let system \eqref{FulSys} satisfy \eqref{FG}, \eqref{pFG}, \eqref{Sform}, \eqref{nres}, and $(n+d,p)\in \mathcal Q_2$ be integers such that assumptions \eqref{as1}, \eqref{asnonlin}, \eqref{Pmas} hold with $\mathcal P'(z_m)<0$. Then for all $\varepsilon>0$ there exist $\delta_0>0$ and $t_0\geq t_\ast$ such that the solution $r(t)$, $\varphi(t)$ of system \eqref{FulSys} with initial data $|r(t_0)-\varrho_m(t_0)|\leq \delta_0$ and $\varphi(t_0)\in \mathbb R$ satisfies 
\begin{gather*} 
\sup_{t\geq t_0}\left\{t^{\vartheta_m}|r(t)-\varrho_m(t)|\right\}<\varepsilon,\quad
\varphi(t)\sim \phi_m(t), \quad  t\to\infty.
\end{gather*}
\end{Th}

The proofs of Lemmas~\ref{Lem1}, \ref{Cor1}, \ref{Lem1u} and Theorems~\ref{Th2}, \ref{LemM}, \ref{ThDop} and~\ref{ThM} are contained in Section~\ref{SQ1}.

\subsection{The case $(n,p)\in\mathcal Q_2$}
In this case, the global behaviour of solutions is determined by the zeros of the coefficient
$\Lambda_p(\rho)$. Consider the following assumption:
\begin{gather}\label{asQ2}
\exists\, \rho_0\in\mathbb R: \quad 0<|\rho_0|< r_\ast, \quad \Lambda_p(\rho_0)=0, \quad \Lambda_p'(\rho_0)\neq 0.
\end{gather}
Then, we have the following:

\begin{Lem}\label{Lem2}
Let $(n,p)\in\mathcal Q_2$ and assumptions \eqref{as1}, \eqref{asQ2} hold with $\Lambda_p'(\rho_0)<0$. Then system \eqref{trSys}, \eqref{trSys2} has a solution $\varrho_2(t)$, $\phi_2(t)$ such that
\begin{gather}\label{sol22}
\varrho_2(t)\sim  \rho_0+
\sum_{k=1}^{\infty} t^{-\frac{k}{q}} \rho_k(\log t), \quad \phi_2(t)   \sim \omega(\rho_0) t, \quad t\to\infty,
\end{gather}
where $\rho_k(\tau)$ are some polynomials in $\tau$. 
\end{Lem}

As in the previous case, the solution $\varrho_2(t)$ is asymptotically stable. 

\begin{Lem}\label{Cor2}
Let $(n,p)\in\mathcal Q_2$ and assumptions \eqref{as1}, \eqref{asQ2} hold with $\Lambda_p'(\rho_0)<0$. Then the solution $\varrho_2(t)$ of equation \eqref{trSys} is exponentially {\rm (}polynomially{\rm)} stable if $p<q$  $(p=q)$. 
\end{Lem}

Thus, we have
\begin{Th}\label{Th3}
Let system \eqref{FulSys} satisfy \eqref{FG}, \eqref{pFG}, \eqref{Sform}, \eqref{nres}, and $(n,p)\in \mathcal Q_2$ be integers such that assumptions \eqref{as1}, \eqref{asQ2} hold with $\Lambda_{p}'(\rho_0)<0$. Then for all $\varepsilon>0$ and $\varkappa\in (0,|\Lambda_p'(\rho_0)|)$ there exist $\delta_0>0$ and $t_0\geq t_\ast$ such that the solution $r(t)$, $\varphi(t)$ of system \eqref{FulSys} with initial data $|r(t_0)-\varrho_2(t_0)|\leq \delta_0$ and $\varphi(t_0)\in \mathbb R$ satisfies 
\begin{gather}\label{rineq2}
\sup_{t\geq t_0}\left\{ t^{\varkappa} |r(t)-\varrho_2(t)|\right\} <\varepsilon,\\
\label{varphiineq2} 
\varphi(t)\sim \phi_2(t), \quad  t\to\infty.
\end{gather}
\end{Th}

Note that, if $\Lambda_p'(\rho_0)>0$, the solution is unstable. Let $\varrho_{2,M}(t)$ be a $M$th partial sum of the series \eqref{sol22}. We have the following:
\begin{Th}\label{Lem2u}
Let system \eqref{FulSys} satisfy \eqref{FG}, \eqref{pFG}, \eqref{Sform}, \eqref{nres}, and $(n,p)\in \mathcal Q_2$ be integers such that assumptions \eqref{as1}, \eqref{asQ2} hold with $\Lambda_p'(\rho_0)>0$. 
Then for all $M\in\mathbb Z_+$ there exists $\varepsilon>0$ such that for all $\delta_0\in (0,\varepsilon)$ the solution of system \eqref{FulSys} with initial data $|r(t_s)-\varrho_{2,M}(t_s)|\leq \delta_0$, $\varphi(t_s)\in\mathbb R$ at some $t_s\geq t_\ast$  satisfies the inequality $t_e^{\nu_0}|r(t_e)-\varrho_{2,M}(t_e)|\geq \varepsilon$ at some $t_e>t_s$.
\end{Th}

Consider the case when, instead of \eqref{asQ2}, the following assumption holds: 
\begin{gather}\label{asQ2n}
\Lambda_p(\rho)\neq 0 \quad \forall |\rho|<r_\ast.
\end{gather}
Then, the solutions of system \eqref{FulSys} exit from a neighbourhood of the equilibrium of the limiting system.
\begin{Th}\label{ThUnst}
Let system \eqref{FulSys} satisfy \eqref{FG}, \eqref{pFG}, \eqref{Sform}, \eqref{nres}, and $(n,p)\in \mathcal Q_2$ be integers such that assumptions \eqref{as1}, \eqref{asQ2n} hold. Then there exist $\varepsilon\in (0,r_\ast)$ such that for all $\delta_0>0$ the solution $r(t)$, $\varphi(t)$ of system \eqref{FulSys} with initial data $|r(t_s)|\leq \delta_0$ and $\varphi(t_s)\in \mathbb R$  at some $t_s\geq t_\ast$ satisfies $|r(t_e)|\geq \varepsilon$ at some $t_e>t_s$.
\end{Th}

The justification of Lemmas~\ref{Lem2}, \ref{Cor2} and Theorem~\ref{Th3}, \ref{Lem2u} and \ref{ThUnst} is contained in Section~\ref{SQ2}.

\subsection{The case $(n,p)\in\mathcal Q_3$}
It can be shown that if $(n,p)\in\mathcal Q_3$, the amplitude of solutions can tend to any $\rho_0\in\mathbb R$ such that
\begin{gather}\label{asQ3}
0<|\rho_0|< r_\ast, \quad \Lambda_n(\rho_0)\neq 0, \quad \Lambda_n'(\rho_0)\neq 0.
\end{gather}
We have
\begin{Lem}\label{Lem3}
Let $(n,p)\in\mathcal Q_3$ and assumption \eqref{as1} hold. Then, for all $\rho_0\in\mathbb R$ satisfying \eqref{asQ3} with $\Lambda_n'(\rho_0)<0$ there exists a solution $\varrho_3(t)$, $\phi_3(t)$ of system \eqref{trSys}, \eqref{trSys2}  such that
\begin{gather}\label{sol31}
\varrho_3(t)\sim \rho_0+ 
\sum_{k=1}^\infty t^{-\frac{n-q+k-1}{q}} \rho_k, \quad \phi_3(t)   \sim \omega(\rho_0) t, \quad t\to\infty,
\end{gather}
where $\rho_k={\hbox{\rm const}}$.
\end{Lem}

In contrast to the previous cases, the solution $\varrho_3(t)$ turns out to be only neutrally stable.

\begin{Lem}\label{Cor3}
Let $(n,p)\in\mathcal Q_3$ and assumptions \eqref{as1}, \eqref{asQ3} hold with $\Lambda_p'(\rho_0)<0$. Then the solution $\varrho_2(t)$ of equation \eqref{trSys} is stable. 
\end{Lem}

Thus, we have the following:
\begin{Th}\label{Th4}
Let system \eqref{FulSys} satisfy \eqref{FG}, \eqref{pFG}, \eqref{Sform}, \eqref{nres}, and $(n,p)\in \mathcal Q_3$ be integers such that assumptions \eqref{as1}, \eqref{asQ3} hold with $\Lambda_{n}'(\rho_0)<0$. Then for all $\varepsilon>0$ there exist $\delta_0>0$ and $t_0\geq t_\ast$ such that the solution $r(t)$, $\varphi(t)$ of system \eqref{FulSys} with initial data $|r(t_0)-\varrho_3(t_0)|\leq \delta_0$ and $\varphi(t_0)\in \mathbb R$ satisfies 
\begin{gather}
\label{rineq3} \sup_{t\geq t_0} |r(t)-\varrho_3(t)| <\varepsilon,\\ 
\label{varphiineq3} \varphi(t)\sim \phi_3(t), \quad  t\to\infty.
\end{gather}
\end{Th}
The proofs of Lemmas~\ref{Lem3}, \ref{Cor3} and Theorem~\ref{Th4} are contained in Section~\ref{SQ3}.

\section{Transformation of variables}\label{Sec3}

To simplify the first equation in system \eqref{FulSys}, consider the transformation of the amplitude $r$ in the following form:
\begin{gather}\label{VN}
V_N(r,\varphi,t)\equiv r+\sum_{k=1}^N t^{-\frac{k}{q}} v_k(r,\varphi,S(t)),
\end{gather}
where the coefficients $v_k(r,\phi,S(t))$ are sought in such a way that the equation for the new variable 
\begin{gather}\label{Rvar}
\rho(t)\equiv V_N(r(t),\varphi(t),S(t),t)
\end{gather} 
does not depend on the ``fast'' angle variables $\varphi$ and $S$ at least in the first $N\in [p,2p)$ terms of its asymptotics as $t\to\infty$. Note that such transformations are usually used in perturbation theory with a small parameter~\cite{BM61,AN84,AKN06}. Differentiating both sides of \eqref{Rvar} with respect to $t$ and taking into account \eqref{FG}, we obtain
\begin{align*}
\frac{d\rho}{dt}\sim & 
\sum_{k=1}^\infty t^{-\frac{k}{q}} 
\left(
\big(\omega(r ) \partial_\varphi +s_0 \partial_S\big) v_k(r ,\varphi,S)+f_k(r ,\varphi,S)+\left(1-\frac{k}{q}\right)v_{k-q}(r ,\varphi,S)\right)
\\
& + \sum_{k=2}^\infty t^{-\frac{k}{q}} \sum_{i+j=k} \left(f_{i}(r ,\varphi,S)\partial_r +r^{-1}g_{i}(r ,\varphi,S)\partial_\varphi+\left(1-\frac{i}{q}\right) s_{i}\partial_S\right) v_{j}(r ,\varphi,S)
\end{align*}
as $t\to\infty$, where it is assumed that $v_i(r ,\varphi,S)\equiv 0$ if $i<1$ or $i>N$, and $s_j=0$ if $j>q$. Comparing asymptotics for $dR/dt$ with the right-hand side of the first equation in \eqref{FulSys3}, we get the following chain of equations:
\begin{gather}\label{veq}
\big(\omega(r ) \partial_\varphi +s_0 \partial_S\big) v_k =\mathcal F_k(r ,\varphi,S), \quad k=1,\dots,N,
\end{gather}
where $\mathcal F_k(r ,\varphi,S)\equiv \Lambda_k(r )-f_k(r ,\varphi,S)+\mathcal R_k(r ,\varphi,S)$. Note that each function $\mathcal R_k(r ,\varphi,S)$ is $2\pi$-periodic with respect to $\varphi$ and $S$, and is defined through $\{v_j(r ,\varphi,S)\}_{j=1}^{k-1}$. For instance,  
\begin{eqnarray*}
\mathcal R_1&\equiv& 0,\\
\mathcal R_2&\equiv &	 v_1 \partial_r  \Lambda_1-\left(1-\frac{2}{q}\right) v_{2-q}-\left(f_1 \partial_r   + r^{-1} g_1 \partial_\varphi  + \left(1-\frac{1}{q}\right) s_1 \partial_S \right)v_1,\\
\mathcal R_3&\equiv &	 \sum_{i+j=3}v_{i} \partial_r  \Lambda_j+\frac{v_1^2}{2}\partial_r ^2 \Lambda_1-
\left(1-\frac{3}{q}\right) v_{3-q}-\sum_{i+j=3}\left(f_i \partial_r   + r^{-1} g_i \partial_\varphi  + \left(1-\frac{1}{q}\right) s_i \partial_S \right)v_j,\\
\mathcal R_k&\equiv &	 \sum_{m_1+2m_2+\dots+i m_i+j=k} C_{i,j,m_1,\dots,m_i}v_{1}^{m_1}\cdots v_i^{m_i} \partial_r ^{m_1+\dots+m_i} \Lambda_j \\
&&-
\left(1-\frac{k}{q}\right) v_{k-q}-\sum_{i+j=k}\left(f_i \partial_r+ r^{-1} g_i \partial_\varphi  + \left(1-\frac{1}{q}\right) s_i \partial_S \right)v_j,
\end{eqnarray*}
where $C_{i,j,m_1,\dots,m_i}={\hbox{\rm const}}$.  Define 
\begin{gather*}
\Lambda_k(r )\equiv \langle f_k(r ,\varphi,S)-\mathcal R_k(r ,\varphi,S)\rangle_{\varphi,S},
\end{gather*}
where the angle brackets denote averaging over $\varphi$ and $S$:
\begin{gather*}
\langle Z(\varphi,S)\rangle_{\varphi,S}:=\frac{1}{(2\pi)^2}\int\limits_0^{2\pi}\int\limits_0^{2\pi} Z(\varphi,S)\,d\varphi\,dS.
\end{gather*} 
In this case, each function $\mathcal F_k(r ,\varphi,S)$ is $2\pi$-periodic with respect to $\varphi$ and $S$ with zero mean. Moreover, $\mathcal F_k(r ,\varphi,S)$ is expanded in Fourier series
\begin{gather*}
\mathcal F_k(r ,\varphi,S)=\sum_{|k_1|+|k_2|\neq 0} \mathcal F_{k,k_1,k_2}(r ) e^{i (k_1 \varphi+k_2 S)}, \\ 
\mathcal F_{k,k_1,k_2}(r )=\frac{1}{(2\pi)^2} \int\limits_0^{2\pi}\int\limits_0^{2\pi} \mathcal F_k(r ,\varphi,S)e^{-i (k_1 \varphi +k_2 S)}\,d\varphi\,dS. 
\end{gather*}
From the properties of the class of perturbations $f $ and $g$ it follows that the Fourier series of $\mathcal F_k(r ,\varphi,S)$  contain a finite number of modes $e^{ik_2S}$. In other words, the set $\mathcal Z:=\{k_2\in\mathbb Z:\mathcal F_{k,k_1,k_2}(r )\not\equiv 0, \ \ \forall\, k, k_1\in\mathbb Z\}$ is finite.
Define the domain
\begin{gather*}
\mathcal D:=\{|r| \leq r_0: \ \ k_1\omega(r )+k_2s_0\neq 0, \quad \forall \,k_1\in\mathbb Z, \quad k_2\in\mathcal Z: \quad |k_1|+|k_2|\neq 0\}.
\end{gather*} 
It follows from \eqref{nres} that there exists $r_\ast\in (0,r _0)$ such that $[-r_\ast,r_\ast]\subset \mathcal D$. 
Hence, $\mathcal D$ is non-empty, and for all $r \in \mathcal D$ the solution of \eqref{veq} is expanded in the Fourier series
\begin{gather*}
v_k(r ,\varphi,S)=\sum_{\substack{ k_1\in\mathbb Z, k_2\in\mathcal Z \\ |k_1|+|k_2|\neq 0}} \frac{\mathcal F_{k,k_1,k_2}(r )}{i (k_1 \omega(r ) +k_2 s_0)}  e^{i (k_1 \varphi +k_2 S)}.
\end{gather*}
Thus, system \eqref{veq} is solvable in the class of $2\pi$-periodic functions in $\varphi$ and $S$ with zero mean. It can easily be checked that 
\begin{gather*}\begin{split}
& \mathcal R_k(r ,\varphi,S)=
	\begin{cases} 
		\mathcal O(r ), & k\leq p,\\
		\mathcal O(1), & p<k\leq N,
	\end{cases}  \quad
	\Lambda_k(r )=
	\begin{cases} 
		\mathcal O(r ), & k< p,\\
		\mathcal O(1), & p\leq k\leq N,
	\end{cases}  \\
	 & v_k(r ,\varphi,S)=
	\begin{cases} 
		\mathcal O(r ), & k< p,\\
		\mathcal O(1), & p\leq k\leq N,
	\end{cases}
	\end{split}
\end{gather*}
as $r \to 0$ uniformly for all $(\varphi,S)\in\mathbb R^2$. In particular, $\Lambda_p(r )=\langle f_p(0,\varphi,S)\rangle_{\varphi,S}+\mathcal O(r )$ as $r \to 0$.

It follows from \eqref{VN} that for all $\epsilon\in (0,r_\ast)$ there exists $t_\ast\geq t_0$ such that
\begin{gather*}
|V_N(r ,\varphi,t)-r |\leq \epsilon, \quad |\partial_r  V_N(r ,\varphi,t)-1|\leq \epsilon
\end{gather*}
for all $|r| \leq r _\ast$, $\varphi\in\mathbb R$ and $t\geq t_\ast$. Hence, the transformation $(r ,\varphi,t)\mapsto (\rho,\varphi,t)$ is invertible for all  $|\rho|\leq \rho_\ast$, $\varphi\in\mathbb R$ and $t\geq t_\ast$ with $\rho_\ast=r _\ast-\epsilon$. Denote by $r =R(\rho,\varphi,t)$ the corresponding inverse transformation. Hence,
\begin{align*}
&\tilde \Lambda_N(\rho,\varphi,t)\equiv \Big(f \partial_r  +(\omega+r^{-1} g) \partial_\varphi +\partial_t \Big)V_N \Big|_{r =R(\rho,\varphi,t)}-\Lambda_N(\rho,t),\\
&\tilde \Omega_{N,1}(\rho,\varphi,t)\equiv \omega(R(\rho,\varphi,t))-\omega(\rho)+r^{-1} \sum_{k=1}^{p-1} t^{-\frac{k}{q}} g_k(R(\rho,\varphi,t),\varphi,S(t),t),\\
&\tilde \Omega_{N,p}(\rho,\varphi,t)\equiv r^{-1}g(R(\rho,\varphi,t),\varphi,S(t),t) - r^{-1} \sum_{k=1}^{p-1} t^{-\frac{k}{q}} g_k(R(\rho,\varphi,t),\varphi,S(t),t).
\end{align*}
Combining this with \eqref{pFG}, we obtain \eqref{LNONest}.

Thus, we obtain the proof of Theorem~\ref{Th1} with $\tilde v_N(r ,\varphi,t)\equiv V_N(r ,\varphi,t)-r$, 
\begin{gather*}
\tilde v_{N,1}(r ,\varphi,t)\equiv \sum_{k=1}^{p-1}t^{-\frac{k-1}{q}}v_k(r,\varphi,S(t)), \quad 
\tilde v_{N,p}(r ,\varphi,t)\equiv \sum_{k=p}^{N}t^{-\frac{k-p}{q}}v_k(r,\varphi,S(t)).
\end{gather*}

\section{The case $(n,p)\in\mathcal Q_1$}\label{SQ1}%%%%%%%%%%%%%%%%%%%%%%%%%%%%%%%%%%%%%%%%%%%%%%%%%%%%%%%%%%%%%%%

\subsection{Proof of Lemma~\ref{Lem1}}
Substituting $\varrho(t)=t^{-\nu_0} \xi(t)$ into equation \eqref{trSys} yields
\begin{gather}\label{xieq1}
\frac{d\xi}{dt}=A_N(\xi,t), \quad A_N(\xi,t)\equiv t^{\nu_0} \Lambda_N\left(t^{-\nu_0}\xi,t\right)+\nu_0 t^{-1}\xi.
\end{gather}
It can easily be checked that the function $A(\xi,t)$ has the following asymptotics:
\begin{gather*} 
A_N(\xi,t)\sim \sum_{k=n}^\infty t^{-\frac{k}{q}} a_k(\xi), \quad t\to\infty, \quad a_n(\xi)\equiv \mu_p+ (\lambda_n+\delta_{n,q}\nu_0) \xi.
\end{gather*}
Consider a particular solution of \eqref{xieq1} tending to a root $\xi_0$ of the equation $a_n(\xi_0)=0$. 

If $n<q$, the asymptotics for a solution of equation \eqref{xieq1} can be constructed in the form of power series with constant coefficients:
\begin{gather*}%\label{xias1}
\xi(t)= \xi_0+\sum_{k=1}^\infty t^{-\frac{k}{q}}\xi_k, \quad \xi_0=-\frac{\mu_p}{\lambda_n}, \quad \xi_k={\hbox{\rm const}}. 
\end{gather*}
Substituting these series into \eqref{xieq1} and grouping the terms of the same power of $t$ yield the following chain of equations for the coefficients $\xi_k$: 
\begin{gather}\label{xikeq}
	-\lambda_n\xi_k=b_k(\xi_0,\xi_1,\dots,\xi_{k-1}), \quad k\geq 1,
\end{gather}
where the right-hand sides $b_k$ are certain polynomial functions of $(\xi_1,\dots,\xi_{k-1})$. In particular, $b_1\equiv a_{n+1}(\xi_0)$, $b_2\equiv a_{n+2}(\xi_0)+a_{n+1}'(\xi_0)\xi_1+a_n''(\xi_0)\xi_1^2/2$. Since $\lambda_n\neq 0$, system \eqref{xikeq} is solvable. 

To prove the existence of a solution to equation \eqref{xieq1} with such asymptotic expansion, consider the function
\begin{gather*} 
\xi_M(t)\equiv \sum_{k=0}^M t^{-\frac{k}{q}}\xi_k
\end{gather*}
with some $M\in\mathbb Z_+$. It can easily be checked that
\begin{gather*}
\frac{d\xi_M(t)}{dt}-A_N(\xi_M(t),t)=\mathcal O\left(t^{-\frac{M+n+2}{q}}\right), \quad t\to\infty.
\end{gather*}
Substituting 
$
\xi(t)=\xi_M(t)+ t^{-{M}/{q}} z(t) 
$
into equation \eqref{xieq1}, we obtain
\begin{gather}\label{zeq1}
\frac{dz}{dt}= \mathcal A(t) z + \mathcal B(z,t)+ \mathcal C(t),
\end{gather}
where
\begin{align*}
\mathcal A(t)& \equiv \partial_\xi A_N(\xi_M(t),t)+\frac{M}{q}t^{-1},\\
\mathcal B(z,t)& \equiv  t^{\frac{M}{q}}
\left[ 
A_N (\xi_M(t)+ t^{-\frac{M}{q}} z,t) - A_N (\xi_M(t),t)
\right]    - \partial_\xi A_N(\xi_M(t),t) z,\\
\mathcal C(t) & \equiv t^{\frac{M}{q}} \left[ 
A_N (\xi_M(t),t)-\frac{d\xi_M(t)}{dt}
\right].
\end{align*}
We see that
\begin{gather}\label{ABCas}
\mathcal A(t)=t^{-\frac{n}{q}}\left(\lambda_n + \mathcal O(t^{-\frac{1}{q}})\right), \quad 
\mathcal B(z,t)=\mathcal O(z^2)\mathcal O(t^{-\frac{n+M+1}{q}}), \quad 
\mathcal C(t)=\mathcal O(t^{-\frac{n+2}{q}})
\end{gather}
as $t\to\infty$ uniformly for all $|z|\leq \zeta_\ast$ with some $\zeta_\ast={\hbox{\rm const}}>0$. Note that if $\mathcal C(t)\equiv 0$, equation \eqref{zeq1} has the trivial solution $z(t)\equiv 0$. Hence, the function $\mathcal C(t)$ plays the role of an external perturbation. Let us prove the stability of the trivial solution with respect to this perturbation. Using $U(z)=|z|$ as a Lyapunov function candidate for equation \eqref{zeq1}, we obtain 
\begin{gather}\label{dU}
\frac{dU}{dt}=\mathcal A(t) |z| +   \left(\mathcal B(z,t) +  \mathcal C(t)\right){\hbox{\rm sgn }}z. 
\end{gather}
It follows from \eqref{ABCas} that there exist $t_1\geq t_\ast$ and $\zeta_1\in (0,\zeta_\ast)$ such that 
\begin{gather*}
\mathcal A(t)\leq -t^{-\frac{n}{q}}\frac{|\lambda_n|}{2}, \quad
   |\mathcal B(z,t)|   \leq t^{-\frac{n}{q}} \frac{|\lambda_n|}{4}|z|, \quad 
 \mathcal C(t) \leq  t^{-\frac{n+1}{q}} \mathcal C_0
\end{gather*}
for all $t\geq t_1$ and $|z|\leq \zeta_1$ with some $\mathcal C_0={\hbox{\rm const}}>0$. Moreover, for all $\varepsilon\in (0,\zeta_1)$ there exist
\begin{gather*}
\delta_\varepsilon=\frac{8 \mathcal C_0 t_\varepsilon^{-\frac 1q}}{|\lambda_n|}<\varepsilon, \quad 
t_\varepsilon=
	\max
		\left\{
			t_1, 
			\left(\frac{16 \mathcal C_0}{\varepsilon |\lambda_n|}\right)^{q}
		\right\}
\end{gather*}
such that 
\begin{gather*}
\frac{dU}{dt}\leq t^{-\frac{n}{q}}\left(-\frac{|\lambda_n|}{4}+ t_\varepsilon^{-\frac{1}{q}}\frac{\mathcal C_0}{\delta_\varepsilon}\right) U \leq - t^{-\frac{n}{q}}\frac{|\lambda_n|}{8} U< 0
\end{gather*}
for all $\delta_\varepsilon\leq |z|\leq \varepsilon$ and $t\geq t_\varepsilon$. Combining this with the inequalities
\begin{gather}\label{supinf}
\sup_{|z|\leq \delta_\varepsilon} U(z)\leq \delta_\varepsilon<\varepsilon=\inf_{|z|=\varepsilon} U(z),
\end{gather}
we see that any solution $z(t)$ of equation \eqref{zeq1} with initial data $|z(t_\varepsilon)|\leq \delta_\varepsilon$ cannot exit from the domain $|z|\leq \varepsilon$ as $t\geq t_\varepsilon$. Hence, for all $M\in\mathbb Z_+$ the solutions $\xi_\ast(t)$ of equation \eqref{xieq1} starting in the vicinity of $\xi_0$ satisfy the estimate 
$\xi_\ast(t)=\xi_M(t)+\mathcal O(t^{-\frac{M}{q}})$ as $ t\to\infty$.
Returning to the original variables, we obtain the existence of the solution $\varrho_1(t)=t^{-\nu_0}\xi_\ast(t)$ with asymptotics \eqref{sol11}, where $\rho_k(\tau)\equiv \xi_k$. Note that similar method for justifying asymptotic expansions at infinity was used in \cite{LK15} for the case when the corresponding linearised system has negative eigenvalues.

Consider now the case $n=q$. In this case, the asymptotic solution is constructed in the form (see~\cite{KF13})
\begin{gather}\label{xias2}
\xi(t)= \xi_0+\sum_{k=1}^\infty t^{-\frac{k}{q}}\xi_k(\log t),  \quad \xi_0=-\frac{\mu_p}{\lambda_n+\nu_0},
\end{gather}
where the coefficients $\xi_k(\tau)$ depend polynomially on $\tau$. Substituting \eqref{xias2} into \eqref{xieq1}, we obtain the chain of differential equations
\begin{gather}\label{xikeq2}
\frac{d\xi_k}{d\tau}-\left(\lambda_n+\nu_0+\frac{k}{q}\right)\xi_k=b_k(\xi_0,\xi_1,\dots,\xi_{k-1}), \quad k\geq 1.
\end{gather}
Since $\lambda_n+\nu_0<0$, it follows that system \eqref{xikeq2} is solvable and 
\begin{gather*}
\xi_k(\tau)\equiv 
	\begin{cases}
		\displaystyle -\frac{b_k(\xi_0,\xi_1(\tau),\dots,\xi_{k-1}(\tau))}{\lambda_n+\nu_0+k/q}, & \displaystyle k<q|\lambda_n+\nu_0|,\\
		\displaystyle \int\limits_{\log t_\ast}^\tau b_k(\xi_0,\xi_1(\theta),\dots,\xi_{k-1}(\theta))\,d\theta, & \displaystyle k=q|\lambda_n+\nu_0|,\\
		\displaystyle \int\limits_{\infty}^\tau e^{-(\lambda_n+\nu_0+k/q)(\theta-\tau)}b_k(\xi_0,\xi_1(\theta),\dots,\xi_{k-1}(\theta))\,d\theta, & \displaystyle k>q|\lambda_n+\nu_0|.
	\end{cases}
\end{gather*}
Define the function 
\begin{gather*}
\xi_M(t)\equiv \xi_0+\sum_{k=1}^M t^{-\frac{k}{q}}\xi_k(\log t), \quad M\in\mathbb Z_+.
\end{gather*}
 It is readily seen that in this case
\begin{gather*}
\frac{d\xi_M(t)}{dt}-A_N(\xi_M(t),t)=\mathcal O\left(t^{-1-\frac{M+1}{q}}\right), \quad t\to\infty.
\end{gather*}
Substituting 
$
\xi(t)=\xi_M(t)+ z(t) 
$
into \eqref{xieq1}, we get \eqref{zeq1} but with the functions 
\begin{align*}
{\mathcal A}(t)& \equiv \partial_\xi A_N(\xi_M(t),t),\\
{\mathcal B}(z,t)& \equiv   
A_N (\xi_M(t)+  z,t) - A_N (\xi_M(t),t)
 - \partial_\xi A_N(\xi_M(t),t) z,\\
{\mathcal C}(t) & \equiv   
A_N (\xi_M(t),t)-\frac{d\xi_M(t)}{dt}
\end{align*}
that satisfy the estimates
\begin{gather}\label{ABCasLem1}
{\mathcal A}(t)=t^{-1}\left(- |\lambda_n+\nu_0| + \mathcal O(t^{-\frac{1}{q}})\right), \quad 
{\mathcal B}(z,t)=\mathcal O(z^2)\mathcal O(t^{-1-\frac{1}{q}}), \quad 
{\mathcal C}(t)=\mathcal O(t^{-1-\frac{M+1}{q}})
\end{gather}
as $t\to\infty$ uniformly for all $|z|\leq \zeta_\ast$ with some $\zeta_\ast={\hbox{\rm const}}>0$. 
As in the previous case, we take $U(z)=|z|$ as a Lyapunov function candidate. In this case, we see that for all $\epsilon\in (0,1)$ there exist $t_1\geq  t_\ast$ and $\zeta_1\in (0,\zeta_\ast)$ such that 
\begin{gather}\label{asest2}
{\mathcal A}(t)\leq -t^{-1}\frac{ |\lambda_n+\nu_0|}{2}, \quad
  |{\mathcal B}(z,t)| \leq t^{-1} \frac{ |\lambda_n+\nu_0|}{4}|z|, \quad 
{\mathcal C}(t) \leq t^{-1-\frac{M+1}{q}}\mathcal C_0 
\end{gather}
for all $t\geq t_1$ and $|z|\leq \zeta_1$ with some $\mathcal C_0={\hbox{\rm const}}>0$.
Hence, for all $\varepsilon\in (0,\zeta_1)$ there exist
\begin{gather*}
\delta_\varepsilon= \frac{8 \mathcal C_0  t_\varepsilon^{-\frac{M+1}{q}} }{ |\lambda_n+\nu_0|}<\varepsilon, \quad 
t_\varepsilon=
	\max
		\left\{
			t_1, 
			\left(\frac{16 \mathcal C_0}{\varepsilon |\lambda_n+\nu_0|}\right)^{\frac{q}{M+1}}
		\right\}
\end{gather*}
such that 
\begin{gather*}
\frac{dU}{dt}\leq t^{-1}\left(-\frac{|\lambda_n+\nu_0|}{4}+ t_\varepsilon^{-\frac{M+1}{q}}\frac{\mathcal C_0}{\delta_\varepsilon}\right) U \leq - t^{-1}\frac{ |\lambda_n+\nu_0|}{8} U< 0
\end{gather*}
for all $\delta_\varepsilon\leq |z|\leq \varepsilon$ and $t\geq t_\varepsilon$. Combining this with \eqref{supinf}, we see that any solution $z(t)$ of equation \eqref{zeq1} with initial data $|z(t_\varepsilon)|\leq \delta_\varepsilon$ cannot exit from the domain $|z|\leq \varepsilon$ as $t\geq t_\varepsilon$. 

From \eqref{dU} and \eqref{asest2} it also follows that 
\begin{gather*}
\frac{dU}{dt}\leq \varepsilon \mathcal C_0 t^{-1-\frac{M+1}{q}}
\end{gather*}
for all $|z|\leq \varepsilon$ and $t\geq t_\varepsilon$. Integrating the last inequality yields $z(t)=\mathcal O(t^{-(M+1)/q})$ as $t\to\infty$. Thus, there exist the solution of equation \eqref{trSys} with asymptotics \eqref{sol11}, where $\rho_k(\tau)\equiv \xi_k(\tau)$.

From equation \eqref{trSys2} it follows that $\phi_1(t)\equiv \int \omega(\varrho_1(t))\,dt$. By using expansion for $\varrho_1(t)$, we obtain asymptotics for $\phi_1(t)$  as $t\to\infty$.

%%%%%%%%%%%%%%%%%%%%%%%%%%%%%%%%%%%%%%%%%%%%%%%%%%%%%%%%%%%%%%%
\subsection{Proof of Lemma~\ref{Cor1}}
Substituting $\varrho(t)=\varrho_{1}(t)+t^{-\nu_0}z$ into equation  \eqref{trSys} yields \eqref{zeq1} with
\begin{align*}
  {\mathcal A}(t)  \equiv \partial_\rho \Lambda_N(\varrho_1(t),t)+\nu_0 t^{-1}, \quad 
  {\mathcal B}(z,t)   \equiv 
t^{\nu_0}\left(\Lambda_N (\varrho_1(t)+t^{-\nu_0}z,t) - \Lambda_N (\varrho_1(t),t)\right)
 - \partial_\rho \Lambda_N(\varrho_1(t),t) z,
\end{align*}
and ${\mathcal C}(t)\equiv 0$. 
It can easily be checked that 
\begin{align*}
  {\mathcal A}(t) =
t^{-\frac{n}{q}}\left(-|\lambda_n+\delta_{n,q}\nu_0| + \mathcal O(t^{-\frac{1}{q}})\right),\quad
  {\mathcal B}(z,t)=\mathcal O(z^2)\mathcal O(t^{-\frac{n+1}{q}})
\end{align*}
as $t\to\infty$ uniformly for all $|\zeta|\leq \zeta_\ast$. Taking $U(z)=|z|$ as a Lyapunov function candidate, it can be seen that there exists $t_1\geq t_\ast$ and $\zeta_1\leq \zeta_\ast$ such that
\begin{gather*}
\frac{dU}{dt}\leq -t^{-\frac{n}{q}}\frac{|\lambda_{n}+\delta_{n,q}\nu_0|}{2} U \leq 0
\end{gather*} 
for all $t\geq t_1$ and $|z|\leq \zeta_1$. Integrating the last inequality, we get
\begin{align*}
&|z(t)|\leq |z(t_1)| \left(\frac{t}{t_1}\right)^{-\frac{|\lambda_n+\nu_0|}{2}}, &\quad & n=q,\\
&|z(t)|\leq |z(t_1)| \exp\left(-\frac{q|\lambda_n|}{2(q-n)} (t^{1-\frac{n}{q}}-t_1^{1-\frac{n}{q}})\right), &\quad & n<q
\end{align*}
as $t\geq t_1$. Returning to the original variables, we see that the solution $\varrho_1(t)$ is polynomially stable if $n=q$ and exponentially stable if $n<q$.

%%%%%%%%%%%%%%%%%%%%%%%%%%%%%%%%%%%%%%%%%%%%%%%%%%%%%%%%%%%%%%%

\subsection{Proof of Theorem~\ref{Th2}}
It follows from Lemma~\ref{Lem1} that if $\tilde \Lambda_N(\rho,\varphi,t)\equiv 0$, the first equation in \eqref{FulSys3} has the solution $\varrho_1(t)$. Let us show that this solution is stable under the perturbation $\tilde \Lambda_N(\rho,\varphi,t)$ uniformly for all $\varphi\in\mathbb R$. Substituting 
\begin{gather}\label{rhosubs1}
\rho(t)=\varrho_1(t)+t^{-\nu_0} z(t)
\end{gather}
into \eqref{FulSys3}, we obtain
\begin{gather}\label{zphisys}
\begin{split}
&\frac{dz}{dt}=\mathcal A(t)z+\mathcal B(z,t)+\mathcal C(z,\varphi,t),\quad 
 \frac{d\varphi}{dt}=\omega(0)+\mathcal D(z,\varphi,t), 
\end{split}
\end{gather}
where
\begin{align*}
\mathcal A(t)& \equiv \partial_\rho \Lambda_N(\varrho_1(t),t)+\nu_0 t^{-1},\\
\mathcal B(z,t)&\equiv t^{\nu_0}\left[\Lambda_N(\varrho_1(t)+t^{-\nu_0} z,t) -\Lambda_N(\varrho_1(t),t)\right]-\partial_\rho \Lambda_N(\varrho_1(t),t)z,\\
\mathcal C(z,\varphi,t)&\equiv t^{\nu_0} \tilde \Lambda_N(\varrho_1(t)+t^{-\nu_0} z,\varphi,t),\\
\mathcal D(z,\varphi,t)& \equiv \omega(\varrho_1(t)+t^{-\nu_0} z)-\omega(0)+\tilde \Omega_N(\varrho_1(t)+t^{-\nu_0} z,\varphi,t).
\end{align*}
It follows from  \eqref{LNONest}, \eqref{as1}, and \eqref{sol11} that 
\begin{align*}
\mathcal A(t)& =t^{-\frac{n}{q}}\left(\lambda_n+\delta_{n,q}\nu_0+\mathcal O(t^{-\frac{1}{q}})\right), \quad 
& \mathcal B(z,t)&=\mathcal O(z^2) \mathcal O(t^{-\frac{n}{q}-\nu_0}), \\
\mathcal C(z,\varphi,t)  & =\mathcal O(t^{\nu_0-\frac{N+1}{q}}), \quad
& \mathcal D(z,\varphi,t)& =\mathcal O(t^{-\frac{1}{q}})
\end{align*}
as $t\to\infty$ uniformly for all $|z|\leq \zeta_\ast$ and $\varphi\in\mathbb R$ with some $\zeta_\ast={\hbox{\rm const}}>0$. Consider $U(z)=|z|$ as a Lyapunov function candidate for system \eqref{zphisys}. The derivative of $U(z)$ along the trajectories of the system is given by 
$
{dU}/{dt}=\mathcal A(t)|z|+ (\mathcal B(z,t)+\mathcal C(z,\varphi,t)) {\hbox{\rm sgn}}z.
$
Taking $N>p+1$, ensures that there exist $t_1\geq t_\ast$ and $\zeta_1\in (0,\zeta_\ast)$ such that
\begin{align*}
 \mathcal A(t)&\leq -t^{-\frac{n}{q}}\frac{|\lambda_n+\delta_{n,q}\nu_0|}{2}, \quad 
&|\mathcal B(z,t) |&\leq t^{-\frac{n}{q}}\frac{|\lambda_n+\delta_{n,q}\nu_0|}{4} |z|, \\
 |\mathcal C(z,\varphi,t)|&\leq  t^{-\frac{n+1}{q}} \mathcal C_0, \quad 
&|\mathcal D(z,\varphi,t)|&\leq  t^{-\frac{1}{q}} \mathcal D_0
\end{align*} 
for all $t\geq t_1$, $|z|\leq \zeta_1$ and $\varphi\in\mathbb R$ with $\mathcal C_0,\mathcal D_0={\hbox{\rm const}}>0$.
Moreover, for all $\varepsilon\in (0,\zeta_1)$ there exist
\begin{gather*}
\delta_\varepsilon=\frac{8\mathcal C_0 t_\varepsilon^{- {1}/{q}}}{ |\lambda_n+\delta_{n,q}\nu_0|}<\varepsilon, \quad
t_\varepsilon=\max\left\{t_1, \left(\frac{16 \mathcal C_0}{\varepsilon |\lambda_{n,q}+\delta_{n,q}\nu_0|}\right)^{q}\right\}
\end{gather*}
such that
\begin{gather*}
\frac{dU}{dt}\leq t^{-\frac{n}{q}}\left(-\frac{|\lambda_n+\delta_{n,q}\nu_0|}{4}+t_\varepsilon^{-\frac{1}{q}}\frac{\mathcal C_0}{\delta_\varepsilon} \right)U\leq -t^{-\frac{n}{q}}  \frac{|\lambda_n+\delta_{n,q}\nu_0|}{8}U< 0
\end{gather*}
for all $\delta_\varepsilon\leq |z|\leq \varepsilon$, $\varphi\in\mathbb R$ and $t\geq t_\varepsilon$. Therefore, any solution $(z(t),\varphi(t))$ of \eqref{zphisys} with initial data $|z(t_\varepsilon)|\leq \delta_\varepsilon$, $\varphi(t_\varepsilon)\in\mathbb R$ satisfies $|z(t)|\leq \varepsilon$ for all $t\geq t_\varepsilon$. Combining this with \eqref{rhosubs0} and \eqref{rhosubs1}, we obtain \eqref{rineq1}.

From the second equation in \eqref{zphisys} it follows that
\begin{align*}
&|\varphi(t)-\varphi(t_\varepsilon)-\omega(0) (t-t_\varepsilon)|\leq \frac{q\mathcal D_0}{q-1} \left(t^{1-\frac{1}{q}}-t_\varepsilon^{1-\frac{1}{q}}\right) \quad & \text{if} &\quad q\neq 1,\\
&|\varphi(t)-\varphi(t_\varepsilon)-\omega(0) (t-t_\varepsilon)|\leq \mathcal D_0 \left(\log {t}-\log {t_\varepsilon}\right) \quad & \text{if}& \quad q= 1
\end{align*}
as $t\geq t_\varepsilon$. Thus, we obtain \eqref{varphiineq1}.

%%%%%%%%%%%%%%%%%%%%%%%%%%%%%%%%%%%%%%%%%%%%%%%%%%%%%%%%%%%%%%%
\subsection{Proof of Theorem~\ref{Lem1u}}
Substituting $\rho(t)=\varrho_{1,M}(t)+t^{-\nu_0}z$ into \eqref{FulSys3}, we obtain \eqref{zphisys} but with 
\begin{align*}
\mathcal A(t)& \equiv \partial_\rho \Lambda_N(\varrho_{1,M}(t),t)+\nu_0 t^{-1},\\
\mathcal B(z,t)&\equiv t^{\nu_0}\left[\Lambda_N(\varrho_{1,M}(t)+t^{-\nu_0} z,t) -\Lambda_N(\varrho_{1,M}(t),t)\right]-\partial_\rho \Lambda_N(\varrho_{1,M}(t),t)z,\\
\mathcal C(z,\varphi,t)&\equiv t^{\nu_0} \tilde \Lambda_N(\varrho_{1,M}(t)+t^{-\nu_0} z,\varphi,t)+
t^{\nu_0}\left[\Lambda_N(\varrho_{1,M}(t),t) -\varrho_{1,M}'(t)\right],\\
\mathcal D(z,\varphi,t)& \equiv \omega(\varrho_{1,M}(t)+t^{-\nu_0} z)-\omega(0)+\tilde \Omega_N(\varrho_{1,M}(t)+t^{-\nu_0} z,\varphi,t)
\end{align*}
such that
\begin{gather*}
{\mathcal A}(t)=t^{-\frac{n}{q}}\left(\lambda_n+\delta_{n,q}\nu_0 + \mathcal O(t^{-\frac{1}{q}})\right), \quad 
{\mathcal B}(z,t)=\mathcal O(z^2)\mathcal O(t^{-\frac{n+1}{q}}), \\ 
{\mathcal C}(z,\varphi,t)=\mathcal O(t^{-\frac{n+1+N-p}{q}})+\mathcal O(t^{-\frac{n+1+M}{q}})
\end{gather*}
 as $t\to\infty$ uniformly for all $|\zeta|\leq \zeta_\ast$ and $\varphi\in\mathbb R$. Consider $U(z)=|z|$ as a Lyapunov function candidate. From \eqref{ABCas} and \eqref{dU} it follows that there exist $t_1\geq t_\ast$ and $\epsilon\in (0,\zeta_\ast)$ such that
\begin{gather*}
\frac{dU}{dt}\geq t^{-\frac{n}{q}}\left(\frac{\lambda_n+\delta_{n,q}\nu_0}{2}U-t^{-\frac{1}{q}} \mathcal C_0  \right)
\end{gather*}
for all $ |z|\leq \epsilon$, $\varphi\in\mathbb R$ and $t\geq t_1$. Hence, for all $\delta\in (0,\epsilon)$ there exists 
\begin{gather*}
t_s=\max\left\{t_1, \left(\frac{4\mathcal C_0}{\delta (\lambda_n+\delta_{n,q}\nu_0)}\right)^{q}\right\}
\end{gather*} 
such that $dU/dt\geq t^{-n/q} (\lambda_n+\delta_{n,q}\nu_0) U/4>0$ for all $\delta\leq|z|\leq \epsilon$, $\varphi\in\mathbb R$ and $t\geq t_s$. Integrating the last inequality, we get $|z(t)|\geq |z(t_s)| (t/t_s)^{(\lambda_n+\nu_0)/4}$ as $t\geq t_s$ if $n=q$. Taking $|z(t_s)|=\delta$, we see that there exists $t_e=t_s (2\epsilon/\delta)^{4/(\lambda_n+\nu_0)}$ such that $|z(t_e)|\geq \epsilon$. Note that similar estimate holds in the case $n<q$. 
 %%%%%%%%%%%%%%%%%%%%%%%%%%%%%%%%%%%%%%%%%%%%%%%%%%%%%%%%%%%%%%%
\subsection{Proof of Theorem~\ref{ThDop}}
From the first equation in \eqref{FulSys3} and assumption \eqref{aslin} it follows that 
\begin{gather*}
\frac{d|\rho|}{dt}=t^{-1}|\rho| \left(\lambda_n+\mathcal O(|\rho|)+\mathcal O(t^{-\frac{1}{q}})\right)+\mathcal O(t^{-\frac{p}{q}})
\end{gather*}
as $\rho\to 0$ and $t\to\infty$ uniformly for all $\varphi\in\mathbb R$. Hence, there exist $\mathcal C_0>0$, $\rho_1\in (0,\rho_\ast)$ and $t_1\geq t_\ast$ such that
\begin{gather*}
\frac{d|\rho|}{dt}\leq t^{-1}\left(-\frac{|\lambda_n|}{2} |\rho|+\mathcal C_0 t^{-\frac{p-q}{q}}\right)
\end{gather*}
for all $|\rho|\leq \rho_1$, $\varphi\in\mathbb R$ and $t\geq t_1$. Moreover, for all $\epsilon\in (0,\rho_1)$ there exist
\begin{gather*}
\delta_\epsilon=\frac{4  \mathcal C_0t_\epsilon^{-\frac{p-q}{q}}}{|\lambda_n|}, \quad t_\epsilon=\max\left\{t_1, \left(\frac{8 \mathcal C_0}{\epsilon|\lambda_n|}\right)^{\frac{q}{p-q}}\right\}
\end{gather*}
such that 
$d|\rho|/dt\leq -t^{-1}|\rho||\lambda_n|/8<0$ for all $\delta_\epsilon\leq |\rho|\leq \epsilon$, $\varphi\in\mathbb R$ and $t\geq t_\epsilon$. Therefore, any solution $\rho(t)$, $\varphi(t)$ of system \eqref{FulSys3} with initial data $|\rho(t_\epsilon)|\leq \delta_\epsilon$, $\varphi(t_\epsilon)\in\mathbb R$ satisfies the inequality $|\rho(t)|\leq \epsilon$ for all $t\geq t_\epsilon$. Returning to the variables $(r,\varphi)$ completes the proof.

%%%%%%%%%%%%%%%%%%%%%%%%%%%%%%%%%%%%%%%%%%%%%%%%%%%%%%%%%%%%%%%

\subsection{Proof of Lemma~\ref{LemM}}
Substituting $\varrho(t)=t^{-\vartheta_m} \xi(t)$ into equation \eqref{trSys} yields \eqref{xieq1} with 
\begin{gather*}
 A_N(\xi,t)\equiv t^{\vartheta_m} \Lambda_N\left(t^{-\vartheta_m}\xi,t\right)+\vartheta_m t^{-1}\xi \sim \sum_{k=n_m}^\infty t^{-\frac{k}{qm}} a_k(\xi), \quad t\to\infty,
\end{gather*}
where
\begin{gather*}
a_{n_m}(\xi)\equiv\mathcal P(\xi), \quad 
n_m=\begin{cases} p(m-1)+n, & m<m_\ast,\\
m(n+d), & m\geq m_\ast.
\end{cases}
\end{gather*}
Note that $n_m/(qm)\leq 1$. The asymptotic solution is constructed in the form
\begin{gather}\label{xim} 
\xi(t)= z_m+\sum_{k=1}^\infty t^{-\frac{k}{qm}}\xi_k(\log t), \quad t\to\infty,
\end{gather}
where the coefficients $\xi_k(\tau)$ depend polynomially on $\tau$ if $n_m=qm$ and $\xi_k(\tau)={\hbox{\rm const}}$ if $n_m<qm$. Substituting \eqref{xim}  into \eqref{xieq1}, we obtain the following chain of equations as $k\geq 1$:
\begin{gather}\label{sysxikm}
\begin{split}
  - \mathcal P'(z_m) \xi_k=b_k(z_m,\xi_1,\dots,\xi_{k-1}), & \quad    n_m<qm,\\
 \frac{d\xi_k}{d\tau}-\left(\mathcal P'(z_m)+\frac{k}{qm}\right)\xi_k=b_k(z_m,\xi_1,\dots,\xi_{k-1}), & \quad    n_m=qm,
\end{split}
\end{gather} 
where  $b_k$ are certain polynomials of $(\xi_1,\dots,\xi_{k-1})$. In particular, $b_1\equiv a_{n_m+1}(z_m)$, $b_2\equiv a_{n_m+2}(z_m)+a_{n_m+1}'(z_m)\xi_1+a_{n_m}''(z_m)\xi_1^2/2$. Since $\mathcal P'(z_m)< 0$, system \eqref{sysxikm} is solvable.

The existence of the solution $\xi_m(t)$ with asymptotics \eqref{solm}, where $\rho_k(\tau)\equiv \xi_k(\tau)$, is proved as in Lemma~\ref{Lem1}. By using  estimates for $\varrho_m(t)=t^{-\vartheta_m} \xi_m(t)$ and  equation \eqref{trSys2}, we obtain asymptotics for $\phi_m(t)$  as $t\to\infty$.

\subsection{Proof of Theorem~\ref{ThM}}
The proof is similar to that of Theorem~\ref{Th2}.

\section{The case $(n,p)\in\mathcal Q_2$}\label{SQ2}

\subsection{Proof of Lemma~\ref{Lem2}}
Consider first the case $n=p<q$. A formal solution of equation \eqref{trSys} is sought in the form of a series
\begin{gather}\label{rhoas1}
\varrho(t)= \rho_0+\sum_{k=1}^\infty t^{-\frac{k}{q}}\rho_k, \quad \rho_k={\hbox{\rm const}}.
\end{gather}
Substituting \eqref{rhoas1} into \eqref{trSys} and equating the terms of the same power of $t$, we get a chain of equations 
\begin{gather}\label{rhokeq}
	-\Lambda_p'(\rho_0)\rho_k=b_k(\rho_0,\rho_1,\dots,\rho_{k-1}), \quad k\geq 1,
\end{gather}
where $b_k(\rho_0,\rho_1,\dots,\rho_{k-1})$ are certain polynomials in $(\rho_1,\dots,\rho_{k-1})$. For instance, $b_1\equiv  \Lambda_{p+1}(\rho_0)$, $b_2\equiv  \Lambda_{p+2}(\rho_0)+\Lambda_{p+1}'(\rho_0)\rho_1+\Lambda_{p}''(\rho_0)\rho_1^2/2$. Since $\Lambda_p'(\rho_0)\neq  0$, system \eqref{rhokeq} is solvable. 

Consider the function
\begin{gather*}%\label{anzM}
\varrho_M(t)\equiv \sum_{k=0}^M t^{-\frac{k}{q}}\rho_k, \quad M\in\mathbb Z_+.
\end{gather*}
By construction,
\begin{gather*}
\frac{d\varrho_M(t)}{dt}-\Lambda_N(\varrho_M(t),t)=\mathcal O\left(t^{-\frac{M+p+2}{q}}\right), \quad t\to\infty.
\end{gather*}
Substituting 
$
\varrho(t)=\varrho_M(t)+ t^{-{M}/{q}} z(t) 
$
into \eqref{trSys}, we obtain equation \eqref{zeq1}, where
\begin{align*}
\mathcal A(t)& \equiv \partial_\rho \Lambda_N(\varrho_M(t),t)+\frac{M}{q}t^{-1}=t^{-\frac{p}{q}}\left(\Lambda_p'(\rho_0) + \mathcal O(t^{-\frac{1}{q}})\right),\\
\mathcal B(z,t)& \equiv  t^{\frac{M}{q}}
\left[ 
\Lambda_N (\varrho_M(t)+ t^{-\frac{M}{q}} z,t) - \Lambda_N (\varrho_M(t),t)
\right]    - \partial_\rho \Lambda_N(\varrho_M(t),t) z=\mathcal O(z^2)\mathcal O(t^{-\frac{p+M+1}{q}}),\\
\mathcal C(t) & \equiv t^{\frac{M}{q}} \left[ 
\Lambda_N (\varrho_M(t),t)-\frac{d\varrho_M(t)}{dt}
\right]=\mathcal O(t^{-\frac{p+2}{q}}).
\end{align*}
as $t\to\infty$ uniformly for all $|z|<\infty$. 
Repeating the argument used in Lemma~\ref{Lem1} with $\Lambda_p'(\rho_0)$ instead of $\lambda_n$, it can be shown that there exists the solution $\varrho_2(t)$ of \eqref{trSys} with asymptotics \eqref{sol22} as $ t\to\infty$.

In the case $n=p=q$, the asymptotic solution is constructed in the following form:
\begin{gather}\label{rhoas2}
\varrho(t)= \rho_0+\sum_{k=1}^\infty t^{-\frac{k}{q}}\rho_k(\log t), \quad t\to\infty,
\end{gather}
where $\rho_k(\tau)$ are polynomials in $\tau$. Substituting \eqref{rhoas2} into \eqref{trSys}, we obtain 
\begin{gather}\label{rhokeq2}
\frac{d\rho_k}{d\tau}-\left(\Lambda_p'(\rho_0)+\frac{k}{q}\right)\rho_k=b_k(\rho_0,\rho_1,\dots,\rho_{k-1}), \quad k\geq 1.
\end{gather}
It can easily be checked that the solution of \eqref{rhokeq2} is given by
\begin{gather*}
\rho_k(\tau)\equiv 
	\begin{cases}
		\displaystyle -\frac{b_k(\rho_0,\rho_1(\tau),\dots,\rho_{k-1}(\tau))}{\Lambda_p'(\rho_0)+k/q}, & \displaystyle k<q|\Lambda_p'(\rho_0)|,\\
		\displaystyle \int\limits_{\tau_\ast}^\tau b_k(\rho_0,\rho_1(\theta),\dots,\rho_{k-1}(\theta))\,d\theta, & \displaystyle k=q|\Lambda_p'(\rho_0)|,\\
		\displaystyle \int\limits_{\infty}^\tau e^{-(\Lambda_p'(\rho_0)+k/q)(\theta-\tau)}b_k(\rho_0,\rho_1(\theta),\dots,\rho_{k-1}(\theta))\,d\theta, & \displaystyle k>q|\Lambda_p'(\rho_0)|,
	\end{cases}
\end{gather*}
with $\tau_\ast=\log t_\ast$. 
Define 
\begin{gather*}
\varrho_M(t)\equiv \rho_0+\sum_{k=1}^M t^{-\frac{k}{q}}\rho_k(\log t), \quad M\in\mathbb Z_+.
\end{gather*}
 In this case,
\begin{gather*}
\frac{d\varrho_M(t)}{dt}-\Lambda_N(\varrho_M(t),t)=\mathcal O\left(t^{-1-\frac{M+1}{q}}\right), \quad t\to\infty.
\end{gather*}
Substituting 
$
\varrho(t)=\varrho_M(t)+ z(t)$
into \eqref{trSys} yields \eqref{zeq1},
where
\begin{align*}
 {\mathcal A}(t)& \equiv \partial_\rho \Lambda_N(\varrho_M(t),t)=
t^{-1}\left(-|\Lambda_p'(\rho_0)| + \mathcal O(t^{-\frac{1}{q}})\right),\\
 {\mathcal B}(z,t)& \equiv 
\Lambda_N (\varrho_M(t)+z,t) - \Lambda_N (\varrho_M(t),t)
 - \partial_\rho \Lambda_N(\varrho_M(t),t) z=\mathcal O(z^2)\mathcal O(t^{-1-\frac{1}{q}}),\\
 {\mathcal C}(t) & \equiv
\Lambda_N (\varrho_M(t),t)-\frac{d\varrho_M(t)}{dt}
=\mathcal O(t^{-1-\frac{M+1}{q}}).
\end{align*}
as $t\to\infty$ uniformly for all $|z|<\infty$. The last asymptotic estimates are similar to \eqref{ABCasLem1}. Hence, by repeating the second part of the proof of Lemma~\ref{Lem1}, we justify the existence of the solution $\varrho_2(t)$ of \eqref{trSys} with asympotics \eqref{sol22}.

From \eqref{trSys2} it follows that $\phi_2(t)\equiv \int \omega(\varrho_2(t))\,dt$. Taking into account asymptotics of $\varrho_2(t)$, we obtain estimate \eqref{sol22} for $\phi_2(t)$ as $t\to\infty$.

%%%%%%%%%%%%%%%%%%%%%%%%%%%%%%%%%%%%%%%%%%%%%%%%%%%%%%%%%%%%%%%%%%%%%%%%%%

\subsection{Proof of Lemma~\ref{Cor2}}
Substituting $\varrho(t)=\varrho_{2}(t)+t^{-\varkappa}z$ with $\varkappa\in (0,|\Lambda_p'(\rho_0)|)$ into equation  \eqref{trSys}, we get \eqref{zeq1}, where
\begin{align*}
  {\mathcal A}(t)  \equiv \partial_\rho \Lambda_N(\varrho_2(t),t)+\varkappa t^{-1}, \quad 
  {\mathcal B}(z,t)   \equiv 
t^{\varkappa}\left(\Lambda_N (\varrho_2(t)+t^{-\varkappa}z,t) - \Lambda_N (\varrho_1(t),t)\right)
 - \partial_\rho \Lambda_N(\varrho_2(t),t) z,
\end{align*}
and ${\mathcal C}(t)\equiv 0$. 
Note that
\begin{align*}
  {\mathcal A}(t) =
t^{-\frac{p}{q}}\left(-\mathcal A_n+ \mathcal O(t^{-\frac{1}{q}})\right),\quad
  {\mathcal B}(z,t)=\mathcal O(z^2)\mathcal O(t^{-\frac{p}{q}-\varkappa})
\end{align*}
as $t\to\infty$ uniformly for all $|\zeta|\leq \zeta_\ast$, where $\mathcal A_n=|\Lambda_p'(\rho_0)|-\delta_{n,q}\varkappa>0$. Taking $U(z)=|z|$ as a Lyapunov function candidate, it can be seen that there exists $t_1\geq t_\ast$ and $\zeta_1\leq \zeta_\ast$ such that
\begin{gather*}
\frac{dU}{dt}\leq -t^{-\frac{p}{q}}\frac{\mathcal A_n}{2} U \leq 0
\end{gather*} 
for all $t\geq t_1$ and $|z|\leq \zeta_1$. Integrating the last inequality, we get
\begin{align*}
&|z(t)|\leq |z(t_1)| \left(\frac{t}{t_1}\right)^{-\frac{\mathcal A_n}{2}}, &\quad & p=q,\\
&|z(t)|\leq |z(t_1)| \exp\left(-\frac{q\mathcal A_n}{2(q-p)} (t^{1-\frac{p}{q}}-t_1^{1-\frac{p}{q}})\right), &\quad & p<q
\end{align*}
as $t\geq t_1$. Returning to the original variables, we see that the solution $\varrho_2(t)$ is polynomially stable if $p=q$ and exponentially stable if $p<q$.
%%%%%%%%%%%%%%%%%%%%%%%%%%%%%%%%%%%%%%%%%%%%%%%%%%%%%%%%%%%%%%%%%%%%%%%%%%

\subsection{Proof of Theorem~\ref{Th3}}
Substituting 
\begin{gather}\label{rhosubs2}
\rho(t)=\varrho_2(t)+t^{-\varkappa} z(t)
\end{gather}
with some $\varkappa\in (0,|\Lambda_n'(\rho_0)|)$ into \eqref{FulSys3} yields
\begin{gather}\label{zphisys2}
\frac{dz}{dt}=\mathcal A(t)z+\mathcal B(z,t)+\mathcal C(z,\varphi,t),\quad
\frac{d\varphi}{dt}=\omega(\rho_0)+\mathcal D(z,\varphi,t), 
\end{gather}
where
\begin{align*}
\mathcal A(t)& \equiv \partial_\rho \Lambda_N(\varrho_2(t),t)+\varkappa t^{-1},\\
\mathcal B(z,t)&\equiv t^{\varkappa}\left[\Lambda_N(\varrho_2(t)+t^{-\varkappa} z,t) -\Lambda_N(\varrho_2(t),t)\right]-\partial_\rho \Lambda_N(\varrho_2(t),t)z,\\
\mathcal C(z,\varphi,t)&\equiv t^{\varkappa} \tilde \Lambda_N(\varrho_2(t)+t^{-\varkappa} z,\varphi,t),\\
\mathcal D(z,\varphi,t)& \equiv \omega(\varrho_2(t)+t^{-\varkappa} z)-\omega(\rho_0)+\tilde \Omega_N(\varrho_2(t)+t^{-\varkappa} z,\varphi,t).
\end{align*}
From \eqref{LNONest}, \eqref{as1} and \eqref{sol22} it follows that 
\begin{align*}
\mathcal A(t)& =t^{-\frac{p}{q}}\left(-\mathcal A_n +o(1)\right), \quad 
& \mathcal B(z,t)&=\mathcal O(z^2) \mathcal O(t^{-\frac{p}{q}-\varkappa}), \\
\mathcal C(z,\varphi,t)  & =\mathcal O(t^{\varkappa-\frac{N+1}{q}}), \quad
& \mathcal D(z,\varphi,t)& =\mathcal O(t^{-\frac{\kappa}{q}})+\mathcal O(t^{-\varkappa})
\end{align*}
as $t\to\infty$ uniformly for all $|z|\leq \zeta_\ast$ and $\varphi\in\mathbb R$ with some $\zeta_\ast={\hbox{\rm const}}>0$, $\kappa\in (0,1)$ and $\mathcal A_n=|\Lambda_p'(\rho_0)|-\delta_{n,q}\varkappa>0$. Consider $U(z)=|z|$ as a Lyapunov function candidate for system \eqref{zphisys2}. The derivative of $U(z)$ along the trajectories of the system is given by 
$
{dU}/{dt}=\mathcal A(t)|z|+(\mathcal B(z,t)+\mathcal C(z,\varphi,t)){\hbox{\rm sgn}}z.
$
Taking $N>p+q\varkappa$ ensures that there exist $t_1\geq t_\ast$ and $\zeta_1\in (0,\zeta_\ast)$ such that
\begin{align*}
 \mathcal A(t)\leq -t^{-\frac{p}{q}}\frac{\mathcal A_n}{2}, \quad 
 |\mathcal B(z,t)|\leq t^{-\frac{p}{q}}\frac{\mathcal A_n}{4} |z|, \quad
 |\mathcal C(z,\varphi,t)|\leq t^{-\frac{p+1}{q}} \mathcal C_0, \quad 
 |\mathcal D(z,\varphi,t)|\leq t^{-\frac{\kappa_1}{q}} \mathcal D_0 
\end{align*} 
for all $t\geq t_1$, $|z|\leq \zeta_1$ and $\varphi\in\mathbb R$ with $\mathcal C_0,\mathcal D_0={\hbox{\rm const}}>0$ and $\kappa_1=\min\{\kappa,\varkappa q\}$. Moreover, for all $\varepsilon\in (0,\zeta_1)$ there exist
\begin{gather*}
\delta_\varepsilon=8 \mathcal A_n^{-1} \mathcal C_0 t_\varepsilon^{-\frac{1}{q}}<\varepsilon, \quad
t_\varepsilon=\max\left\{t_1, \left(\frac{16 \mathcal C_0}{\varepsilon \mathcal A_n}\right)^{q}\right\}
\end{gather*}
such that
\begin{gather*}
\frac{dU}{dt}\leq t^{-\frac{p}{q}}\left(-\frac{\mathcal A_n}{4}+\frac{\mathcal C_0}{\delta_\varepsilon} t_\varepsilon^{-\frac{1}{q}}\right)U\leq -t^{-\frac{p}{q}}  \frac{\mathcal A_n}{8}U< 0
\end{gather*}
for all $\delta_\varepsilon\leq |z|\leq \varepsilon$, $\varphi\in\mathbb R$ and $t\geq t_\varepsilon$. Therefore, any solution $(z(t),\varphi(t))$ of \eqref{zphisys} with initial data $|z(t_\varepsilon)|\leq \delta_\varepsilon$, $\varphi(t_\varepsilon)\in\mathbb R$ satisfies $|z(t)|\leq \varepsilon$ for all $t\geq t_\varepsilon$. Combining this with \eqref{rhosubs0} and \eqref{rhosubs2}, we obtain \eqref{rineq2}.

From the second equation of system \eqref{zphisys2} it follows that
\begin{align*}
 |\varphi(t)-\varphi(t_\varepsilon)-\omega(\rho_0) (t-t_\varepsilon)|\leq \frac{q\mathcal D_0}{q-\kappa_1} \left(t^{1-\frac{\kappa_1}{q}}-t_\varepsilon^{1-\frac{\kappa_1}{q}}\right), \quad t\geq t_\varepsilon.
\end{align*}
Thus, we obtain \eqref{varphiineq2}.
 
%%%%%%%%%%%%%%%%%%%%%%%%%%%%%%%%%%%%%%%%%%%%%%%%%%%%%%%%%%%%%%%%%5

\subsection{Proof of Theorem~\ref{Lem2u}}
Substituting 
$
\rho(t)=\varrho_{2,M}(t)+ z(t)$
into \eqref{FulSys3}, we obtain \eqref{zphisys} with
\begin{align*}
\mathcal A(t)& \equiv \partial_\rho \Lambda_N(\varrho_{2,M}(t),t),\\
\mathcal B(z,t)&\equiv \Lambda_N(\varrho_{2,M}(t)+ z,t) -\Lambda_N(\varrho_{2,M}(t),t)-\partial_\rho \Lambda_N(\varrho_{2,M}(t),t)z,\\
\mathcal C(z,\varphi,t)&\equiv \tilde \Lambda_N(\varrho_{2,M}(t)+z,\varphi,t)+
\Lambda_N(\varrho_{2,M}(t),t) -\varrho_{2,M}'(t),\\
\mathcal D(z,\varphi,t)& \equiv \omega(\varrho_{2,M}(t)+ z)-\omega(0)+\tilde \Omega_N(\varrho_{2,M}(t)+z,\varphi,t).
\end{align*}
Note that the following estimates hold:
\begin{gather*}
{\mathcal A}(t)=t^{-\frac{p}{q}}\left( \Lambda_p'(\rho_0)  + \mathcal O(t^{-\frac{1}{q}})\right), \quad 
{\mathcal B}(z,t)=\mathcal O(z^2)\mathcal O(t^{-\frac{p+1}{q}}), \quad 
{\mathcal C}(z,\varphi,t)=\mathcal O(t^{-\frac{N+1}{q}})+\mathcal O(t^{-\frac{p+M+1}{q}})
\end{gather*}
 as $t\to\infty$ uniformly for all $|\zeta|\leq \zeta_\ast$ and $\varphi\in\mathbb R$. Consider $U(z)=|z|$ as a Lyapunov function candidate. Recall that $N\geq p$. Hence, there exist $\mathcal C_0>0$, $t_1\geq t_\ast$ and $\epsilon\in (0,\zeta_\ast)$ such that
\begin{gather*}
\frac{dU}{dt}\geq t^{-\frac{p}{q}}\left(\frac{\Lambda_p'(\rho_0)}{2}U-  t^{-\frac{1}{q}}\mathcal C_0\right)
\end{gather*}
for all $ |z|\leq \epsilon$ and $t\geq t_1$. Thus, for all $\delta\in (0,\epsilon)$ there is 
\begin{gather*}
t_s=\max\left\{t_1, \left(\frac{4\mathcal C_0}{\delta \Lambda_p'(\rho_0)}\right)^{q}\right\}
\end{gather*} 
such that $dU/dt\geq t^{-p/q} \Lambda_p'(\rho_0) U/4>0$ for all $\delta\leq|z|\leq \epsilon$ and $t\geq t_s$. Integrating the last inequality, we get $|z(t)|\geq |z(t_s)| (t/t_s)^{\Lambda_p'(\rho_0)/4}$ as $t\geq t_s$ if $p=q$. Taking $|z(t_s)|=\delta$, we see that there exists $t_e=t_s (2\epsilon/\delta)^{4/\Lambda_p'(\rho_0)}$ such that $|z(t_e)|\geq \epsilon$. The similar estimate holds in the case $p<q$. 

%%%%%%%%%%%%%%%%%%%%%%%%%%%%%%%%%%%%%%%%%%%%%%%%%%%%%%%%%%%%%%%%%%%%%%%%%%
\subsection{Proof of Theorem~\ref{ThUnst}}
 Consider the first equation of system \eqref{FulSys3}. From \eqref{asQ2n} it follows that there exists $c>0$ such that $|\Lambda_p(\rho)|\geq c$ for all $|\rho|<r_\ast$. Hence there exists $t_1\geq t_\ast$ such that 
\begin{gather*}
 \frac{d\rho}{dt}\geq t^{-\frac{p}{q}} \frac{c}{2} \quad \text{if} \quad \Lambda_p(\rho)>0;\quad 
 \frac{d\rho}{dt}\leq - t^{-\frac{p}{q}} \frac{c}{2} \quad \text{if} \quad \Lambda_p(\rho)<0
\end{gather*}
as $t\geq t_1$ for all $|\rho|<r_\ast$ and $\varphi\in\mathbb R$. If $p=q$, then integrating the last inequality yields 
\begin{gather*}
|\rho(t)-\rho(t_1)|\geq \frac{c}{2}\log \frac{t}{t_1}, \quad t\geq t_1.
\end{gather*}
Therefore, there exists $\epsilon \in (0,r_\ast)$ such that for all $\delta\in (0,\epsilon)$ the solution of system \eqref{FulSys3} with initial data $|\rho(t_1)|\leq \delta$ and $\varphi(t_1)\in\mathbb R$ satisfies the inequality $|\rho(t_e)|\geq \epsilon$ at $t_e=t_1 \exp(2(2\epsilon-\delta)/c)$. It can easily be checked that similar inequality holds in the case $p<q$. Returning to the original variables $(r,\varphi)$ completes the proof.

%%%%%%%%%%%%%%%%%%%%%%%%%%%%%%%%%%%%%%%%%%%%%%%%%%%%%%%%%%%%%%%%%%%%%%%%%%
\section{The case $(n,p)\in\mathcal Q_3$}\label{SQ3}

\subsection{Proof of Lemma~\ref{Lem3}}
A formal solution of equation \eqref{trSys} is constructed in the form
\begin{gather}\label{rhoas3}
\varrho(t)= \rho_0+\sum_{k=1}^\infty t^{-\frac{n-q+k-1}{q}}\rho_k, \quad \rho_k={\hbox{\rm const}}.
\end{gather}
Substituting \eqref{rhoas3} into \eqref{trSys}, we obtain a chain of equations for $\rho_k$: 
\begin{gather}\label{rhokeq3}
	-\frac{n-q+k-1}{q}\rho_k=b_k, \quad k\geq 1,
\end{gather}
where $b_k$ are expressed through $\rho_0, \dots, \rho_{k-1}$. For example, $b_1\equiv  \Lambda_{n}(\rho_0)$, $b_2\equiv  \Lambda_{n+1}(\rho_0)+\Lambda_{n}'(\rho_0)\rho_1$, etc. Since $n>q$, we see that system \eqref{rhokeq3} is solvable. Consider the function
\begin{gather*}%\label{anzM3}
\varrho_M(t)\equiv \sum_{k=0}^M t^{-\frac{k}{q}}\rho_k, \quad M\in\mathbb Z_+.
\end{gather*}
It can easily be checked that
\begin{gather*}
\frac{d\varrho_M(t)}{dt}-\Lambda_N(\varrho_M(t),t)=\mathcal O\left(t^{-\frac{M+n+1}{q}}\right), \quad t\to\infty.
\end{gather*}
Substituting 
$
\varrho(t)=\varrho_M(t)+z(t) 
$
into \eqref{trSys}, we obtain equation \eqref{zeq1} with
where
\begin{align*}
\mathcal A(t)& \equiv \partial_\rho \Lambda_N(\varrho_M(t),t)=t^{-\frac{n}{q}}\left(\Lambda_n'(\rho_0) + \mathcal O(t^{-\frac{1}{q}})\right),\\
\mathcal B(z,t)& \equiv   
\Lambda_N (\varrho_M(t)+  z,t) - \Lambda_N (\varrho_M(t),t)   - \partial_\rho \Lambda_N(\varrho_M(t),t) z=\mathcal O(z^2)\mathcal O(t^{-\frac{n+1}{q}}),\\
\mathcal C(t) & \equiv  
\Lambda_N (\varrho_M(t),t)-\frac{d\varrho_M(t)}{dt}
=\mathcal O(t^{-\frac{M+n+1}{q}}).
\end{align*}
as $t\to\infty$ uniformly for all $|z|<\infty$. 
Then, repeating the argument used in Lemma~\ref{Lem1} with $\Lambda_n'(\rho_0)$ instead of $\lambda_n$, it can be shown that there exists the solution $\varrho_3(t)$ of equation \eqref{trSys} with asymptotics \eqref{sol31} as $ t\to\infty$.

%%%%%%%%%%%%%%%%%%%%%%%%%%%%%%%%%%%%%%%%%%%%%%%%%%%%%%%%%%%%%%%%%%%%%%%%%%

\subsection{Proof of Lemma~\ref{Cor3}}
Substituting  $\rho(t)=\varrho_3(t)+z(t)$ into equation \eqref{trSys} yields \eqref{zeq1} with
\begin{gather*}
\mathcal A(t)  \equiv \partial_\rho \Lambda_N(\varrho_3(t),t) ,\quad
\mathcal B(z,t) \equiv  \Lambda_N(\varrho_3(t)+  z,t) -\Lambda_N(\varrho_3(t),t) -\partial_\rho \Lambda_N(\varrho_3(t),t)z
\end{gather*}
and $\mathcal C(t)\equiv 0$. It can easily be checked that
\begin{align*}
  {\mathcal A}(t) =
t^{-\frac{n}{q}}\left(-|\Lambda_n'(\rho_0)|+ \mathcal O(t^{-\frac{1}{q}})\right),\quad
  {\mathcal B}(z,t)=\mathcal O(z^2)\mathcal O(t^{-\frac{p}{q}-\varkappa})
\end{align*}
as $t\to\infty$ uniformly for all $|\zeta|\leq \zeta_\ast$. Taking $U(z)=|z|$ as a Lyapunov function candidate, it can be seen that there exist $t_1\geq t_\ast$ and $\zeta_1\leq \zeta_\ast$ such that
\begin{gather*}
\frac{dU}{dt}\leq -t^{-\frac{n}{q}}\frac{|\Lambda_n'(\rho_0)|}{2} U \leq 0
\end{gather*} 
for all $t\geq t_1$ and $|z|\leq \zeta_1$. It follows that $|z(t)|\leq |z(t_1)|$ as $t\geq t_1$. 
Hence, for all $\varepsilon\in (0,\zeta_1)$ there exists $\delta<\varepsilon$ such that the solution $z(t)$ with initial data $|z(t_1)|\leq \delta$ satisfies the inequality $|z(t)|\leq \varepsilon$ as $t\geq t_1$. Returning to the original variables, we obtain the stability of the solution $\varrho_3(t)$.
 
%%%%%%%%%%%%%%%%%%%%%%%%%%%%%%%%%%%%%%%%%%%%%%%%%%%%%%%%%%%%%%%%%%%%%%%%%%

\subsection{Proof of Theorem~\ref{Th4}}
Substituting  $\rho(t)=\varrho_3(t)+z(t)$ into system \eqref{FulSys3}, we obtain \eqref{zphisys2} with
\begin{align*}
\mathcal A(t)& \equiv \partial_\rho \Lambda_N(\varrho_3(t),t) ,\\
\mathcal B(z,t)&\equiv  \Lambda_N(\varrho_3(t)+  z,t) -\Lambda_N(\varrho_3(t),t) -\partial_\rho \Lambda_N(\varrho_3(t),t)z,\\
\mathcal C(z,\varphi,t)&\equiv \tilde \Lambda_N(\varrho_3(t)+ z,\varphi,t),\\
\mathcal D(z,\varphi,t)& \equiv \omega(\varrho_3(t)+  z)-\omega(\rho_0)+\tilde \Omega_N(\varrho_3(t)+  z,\varphi,t).
\end{align*}
 It follows from \eqref{LNONest}, \eqref{as1} and \eqref{sol31}   that 
\begin{align*}
\mathcal A(t)& =t^{-\frac{n}{q}}\left(- |\Lambda_n'(\rho_0)|  +\mathcal O(t^{-\frac{1}{q}})\right), \quad 
& \mathcal B(z,t)&=\mathcal O(z^2) \mathcal O(t^{-\frac{n}{q}}), \\
\mathcal C(z,\varphi,t)  & =\mathcal O(t^{ -\frac{N+1}{q}}), \quad
& \mathcal D(z,\varphi,t)& =\mathcal O(t^{-\frac{1}{q}})
\end{align*}
as $t\to\infty$ uniformly for all $|z|\leq \zeta_\ast$ and $\varphi\in\mathbb R$ with some $\zeta_\ast={\hbox{\rm const}}>0$.
 
Consider $U(z)=|z|$ as a Lyapunov function candidate. Hence, 
$
{dU}/{dt}=\mathcal A(t)|z|+(\mathcal B(z,t)+\mathcal C(z,\varphi,t)){\hbox{\rm sgn}}z.
$
Taking $N\geq n+1$ ensures that there exist $t_1\geq t_\ast$ and $\zeta_1\in (0,\zeta_\ast)$ such that
\begin{align*}
 \mathcal A(t)&\leq -t^{-\frac{n}{q}}\frac{|\Lambda_n'(\rho_0)|}{2}, \quad 
 |\mathcal B(z,t)|&\leq t^{-\frac{n}{q}}\frac{|\Lambda_n'(\rho_0)|}{4} |z|, \quad
 |\mathcal C(z,\varphi,t)|&\leq  t^{-\frac{n+1}{q}} \mathcal C_0, \quad 
 |\mathcal D(z,\varphi,t)|&\leq t^{-\frac{1}{q}} \mathcal D_0
\end{align*} 
for all $t\geq t_1$, $|z|\leq \zeta_1$ and $\varphi\in\mathbb R$ with $\mathcal C_0,\mathcal D_0={\hbox{\rm const}}>0$.
Therefore, for all $\varepsilon\in (0,\zeta_1)$ there exist
\begin{gather*}
\delta_\varepsilon=\frac{8\mathcal C_0 t_\varepsilon^{- \frac{1}{q}}}{ |\Lambda_n'(\rho_0)|}<\varepsilon, \quad
t_\varepsilon=\max\left\{t_1, \left(\frac{16 \mathcal C_0}{\varepsilon |\Lambda_n'(\rho_0)|}\right)^{q}\right\}
\end{gather*}
such that
\begin{gather*}
\frac{dU}{dt}\leq t^{-\frac{n}{q}}\left(-\frac{|\Lambda_n'(\rho_0)|}{4}+\frac{\mathcal C_0}{\delta_\varepsilon} t_\varepsilon^{-\frac{1}{q}}\right)U\leq -t^{-\frac{n}{q}}  \frac{|\Lambda_n'(\rho_0)|}{8}U< 0
\end{gather*}
for all $\delta_\varepsilon\leq |z|\leq \varepsilon$, $\varphi\in\mathbb R$ and $t\geq t_\varepsilon$. It follows that any solution $(z(t),\varphi(t))$ of system \eqref{zphisys} with initial data $|z(t_\varepsilon)|\leq \delta_\varepsilon$, $\varphi(t_\varepsilon)\in\mathbb R$ satisfies $|z(t)|\leq \varepsilon$ for all $t\geq t_\varepsilon$. Returning to the original variables, we obtain \eqref{rineq3}.

Using the second equation of system \eqref{zphisys2} and the estimate of the function $\mathcal D(z,\varphi,t)$, we get \eqref{varphiineq3}.

\section{Examples}\label{SEx}

\subsection{Example 1} Consider the system 
\begin{gather}\label{Ex1}
\begin{split}
&\frac{dr}{dt}=- t^{-\frac{h}{2}} \left(\alpha(S(t))\cos\varphi-\beta(S(t))\sin\varphi\right) r \sin\varphi +t^{-\frac{p}{2}} \gamma(S(t))\sin^2\varphi, \\
&\frac{d\varphi}{dt}=1-t^{-\frac{h}{2}}\left(\alpha(S(t))\cos\varphi-\beta(S(t))\sin\varphi\right) \cos\varphi +t^{-\frac{p}{2}} \gamma(S(t))\frac{\sin\varphi\cos\varphi}{r},
\end{split}
\end{gather}
where
$\alpha(S)\equiv \alpha_0+\alpha_1 \sin S$, 
$\beta(S)\equiv\beta_0+\beta_1 \sin S$, 
$\gamma(S)\equiv\gamma_0+\gamma_1 \sin S$,
$S(t)\equiv \sqrt 2 t+s_2 \log t$, $\alpha_j, \beta_j, \gamma_j, s_2\in\mathbb R$, $h,p\in\mathbb Z_+$. Note that system \eqref{Ex1} has the form \eqref{FulSys} with $\omega(r)\equiv 1$, 
 \begin{align*}
&f(r,\varphi,t)\equiv - t^{-\frac{h}{2}} \left(\alpha(S(t))\cos\varphi-\beta(S(t))\sin\varphi\right) r \sin\varphi +t^{-\frac{p}{2}} \gamma(S(t))\sin^2\varphi,\\
&g(r,\varphi,t)\equiv -t^{-\frac{h}{2}}\left(\alpha(S(t))\cos\varphi-\beta(S(t))\sin\varphi\right) r\cos\varphi +t^{-\frac{p}{2}} \gamma(S(t))\sin\varphi\cos\varphi,
\end{align*}
and satisfies \eqref{FG}, \eqref{pFG}, \eqref{Sform} and \eqref{nres} with $q=2$, $f_{p,0,0}=\gamma_0/2$, $s_0=\sqrt 2$, $s_1=0$. 
It can easily be checked that system \eqref{Ex1} in the Cartesian coordinates $x=r\cos\varphi$, $y=-r\sin\varphi$ takes the form
\begin{gather*}%\label{Ex1car}
\frac{dx}{dt}=y, \quad \frac{dy}{dt}=-x+t^{-\frac{h}{2}}G_1(x,y,S(t))+t^{-\frac{p}{2}}G_2(x,y,S(t))
\end{gather*}
with
\begin{gather*}
 G_1(x,y,S)\equiv \alpha(S)x+\beta(S)y, \quad G_2(x,y,S)\equiv \frac{\gamma(S) y}{\sqrt{x^2+y^2}}.
\end{gather*}
Note also that system \eqref{Ex0} from Section~\ref{PS} corresponds to system \eqref{Ex1} with 
$\alpha_0=\alpha_1=\beta_1=\gamma_0=s_2=0$, $h=2$ and $\beta_0=\lambda$.

1. Let $h=2$ and $p=3$. Then the transformation described in Theorem~\ref{Th1} and constructed in Section~\ref{Sec3} reduces the system to \eqref{FulSys3} with
\begin{gather*}
\Lambda_1(\rho)\equiv 0, \quad \Lambda_2(\rho)\equiv \frac{\beta_0 \rho}{2}, \quad \Lambda_3(\rho)\equiv \frac{\gamma_0}{2}.
\end{gather*} 
In this case, $n=q=2$, $\mu_p=\gamma_0/2$ and $(n,p)\in\mathcal Q_1$. We see that assumption \eqref{aslin} holds with $\lambda_n=\beta_0/2$. Hence, $\nu_0=1/2$. It follows from Theorem~\ref{Th2} that if $\beta_0<-1$, then there exists a stable regime for solutions of system \eqref{Ex1} with $r(t)\approx  \varrho_1(t)$  and $\varphi(t)\sim t$, where $\varrho_1(t)$ is a polynomially stable solution of the corresponding truncated equation \eqref{trSys} with asymptotics \eqref{sol11}. In particular, $\varrho_1(t)\sim -t^{-1/2}\gamma_0/(\beta_0+1)$ as $t\to\infty$. From Theorem~\ref{Lem1u} it follows that this regime is unstable if $\beta_0>-1$. However, from Theorem~\ref{ThDop} it follows that $r(t)$ remains close to zero if $\beta_0<0$ (see Fig.~\ref{FigEx11}).
\begin{figure}
\centering
\subfigure[]{
 \includegraphics[width=0.3\linewidth]{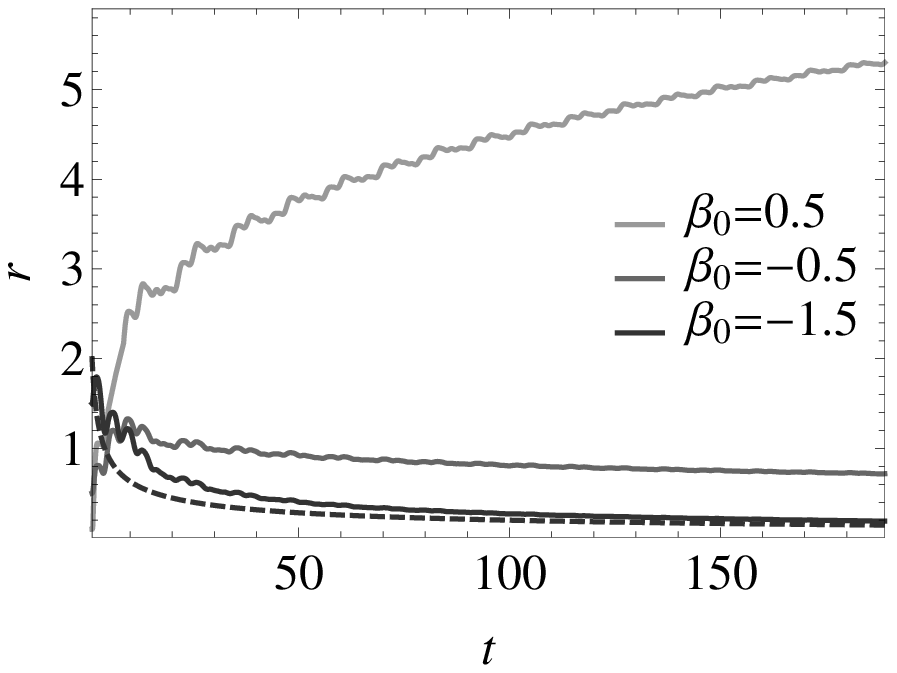}
}
\hspace{1ex}
 \subfigure[]{
 	\includegraphics[width=0.3\linewidth]{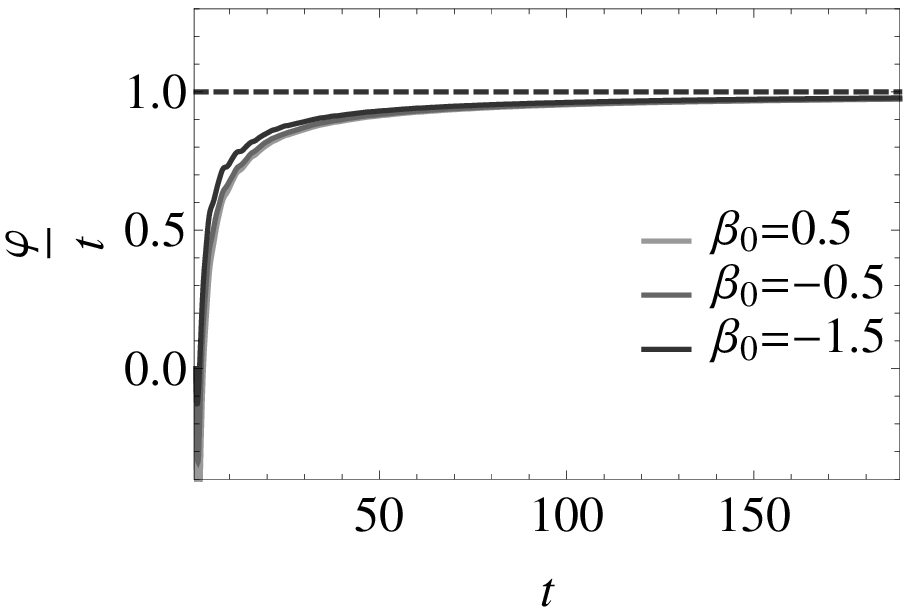}
}
\caption{\small The evolution of $r(t)$ and $\varphi(t)/t$ for solutions to system \eqref{Ex1} with $h=2$, $p=3$, $\alpha_0=\alpha_1=\beta_1=\gamma_0=\gamma_1=s_2=1$ and different values of the parameter $\beta_0$. The dashed curves correspond to (a) $r(t)\equiv 2t^{-1/2}$, (b) $\varphi(t)/t\equiv 1$.} \label{FigEx11}
\end{figure}

2. Now, let $h=p=2$. The change of the variable \eqref{rhosubs0} reduces system \eqref{Ex1} to the form \eqref{FulSys3} 
with
\begin{gather*}
\Lambda_1(\rho)\equiv 0, \quad \Lambda_2(\rho)\equiv \frac{\gamma_0+\beta_0 \rho}{2}.
\end{gather*} 
In this case, $n=q=2$, $(n,p)\in\mathcal Q_2$, and assumption \eqref{asQ2} holds with $\rho_0=-\gamma_0/\beta_0$ and $\Lambda_p'(\rho_0)=\beta_0/2$. It follows from Theorem~\ref{Th3} that if $\beta_0<0$, then there exists a stable regime for solutions of system \eqref{Ex1} such that $r(t)\sim \varrho_2(t)$ and $\varphi(t)\sim t$ as $t\to\infty$, where $\varrho_2(t)$ is a polynomially stable solution of the corresponding truncated equation \eqref{trSys} with asymptotics \eqref{sol22}. We see that $\varrho_2(t)\sim -\gamma_0/\beta_0$ as $t\to\infty$. It follows from Theorem~\ref{Lem2u} that this asymptotic regime becomes unstable if $\beta_0>0$. If $\beta_0=0$ and $\gamma_0\neq 0$, then assumption \eqref{asQ2n} holds. In this case, it follows from Theorem~\ref{ThUnst} that the solutions of system \eqref{Ex1} leave the neighbourhood of zero (see Fig.~\ref{FigEx12}).
\begin{figure}
\centering
\subfigure[]{
 \includegraphics[width=0.3\linewidth]{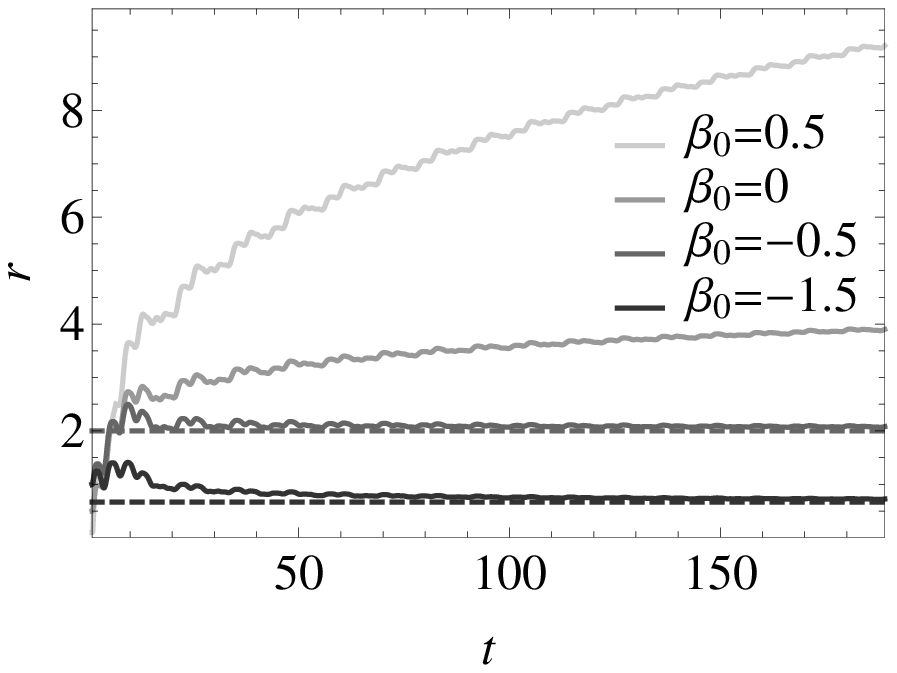}
}
\hspace{1ex}
 \subfigure[]{
 	\includegraphics[width=0.3\linewidth]{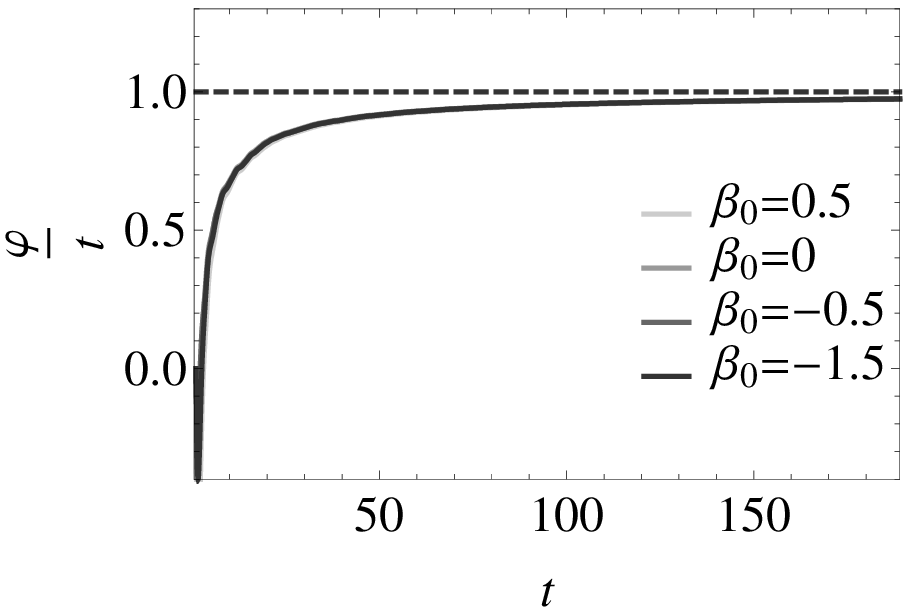}
}
\caption{\small The evolution of $r(t)$ and $\varphi(t)/t$ for solutions to system \eqref{Ex1} with $h=p=2$, $\alpha_0=\alpha_1=\beta_1=\gamma_0=\gamma_1=s_2=1$ and different values of the parameter $\beta_0$. The dashed curves correspond to (a) $r(t)\equiv 2/3$ and $r(t)\equiv 2$, (b) $\varphi(t)/t\equiv 1$.} \label{FigEx12}
\end{figure}

3. Let $h=2$ and $p=1$. Then system \eqref{Ex1} is transformed to \eqref{FulSys3} with $\Lambda_1(\rho)\equiv \gamma_0/2$. In this case, $n=p=1$ and $(n,p)\in\mathcal Q_2$. If $\gamma_0\neq 0$, then assumption \eqref{asQ2n} holds, and it follows from Theorem~\ref{ThUnst} that the solutions of system \eqref{Ex1} leave the neighbourhood of zero.

4. Finally, let $h=p=3$. In this case, it can easily be checked that system \eqref{Ex1} is transformed to \eqref{FulSys3} with 
\begin{gather*}
\Lambda_1(\rho)\equiv \Lambda_2(\rho)=0, \quad \Lambda_3(\rho)\equiv \frac{\gamma_0+\beta_0\rho}{2}.
\end{gather*}
It follows that assumption \eqref{as1} holds with $n=3>q$, and $(n,p)\in\mathcal Q_3$. We see that $\Lambda_3(\rho)\not\equiv 0$ and $\Lambda_3'(\rho)\equiv \beta_0/2$. Hence, assumption \eqref{asQ3} holds if $\rho_0\neq -\gamma_0/\beta_0$. Note that Lemma~\ref{Lem3} ensures that the corresponding truncated equation \eqref{trSys} has a stable solution $\varrho_3(t)\sim \rho_0$ if $\beta_0<0$. From Theorem~\ref{Th4} it follows that there exists a regime for solutions of system \eqref{Ex1} such that $r(t)\approx \varrho_3(t)$ and $\varphi(t)\sim t$ as $t\to\infty$  (see Fig.~\ref{FigEx13}). In this case, the dynamics of the perturbed system \eqref{Ex1} is similar to the dynamics of the corresponding limiting system.
\begin{figure}
\centering{
 \includegraphics[width=0.3\linewidth]{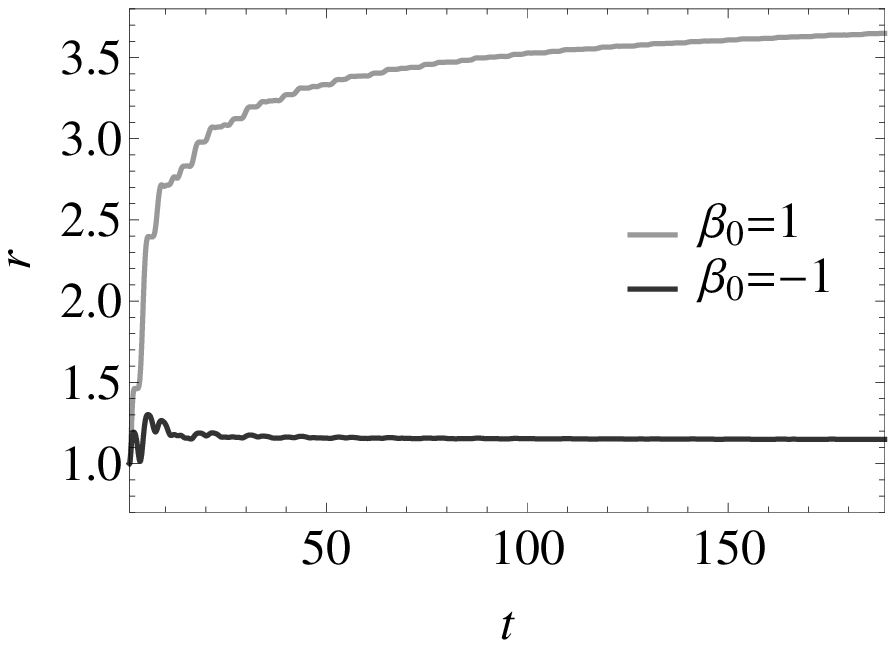}
}
\hspace{1ex}{
 	\includegraphics[width=0.3\linewidth]{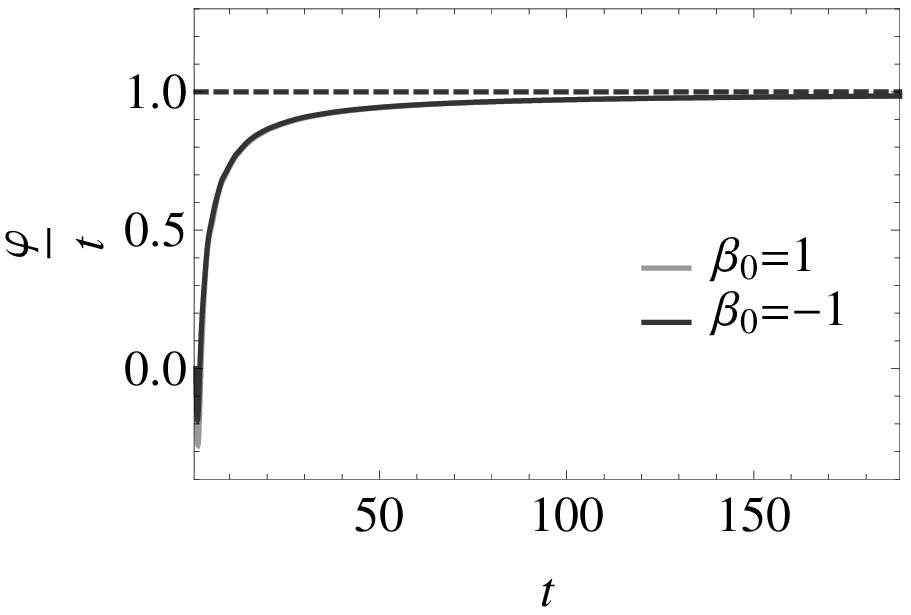}
}
\caption{\small The evolution of $r(t)$ and $\varphi(t)/t$ for solutions to system \eqref{Ex1} with $h=p=3$, $\alpha_0=\alpha_1=\beta_1=\gamma_0=\gamma_1=s_2=1$, with the same initial data and different values of the parameter $\beta_0$. The dashed curve corresponds to $\varphi(t)/t\equiv 1$.} \label{FigEx13}
\end{figure}

\subsection{Example 2} Consider a system in the Cartesian coordinates
\begin{gather}\label{Ex2C}
\frac{dx}{dt}=y, \quad \frac{dy}{dt}=-\sin x + G(x,y,S(t),t),
\end{gather}
where 
\begin{gather*}
G(x,y,S,t)\equiv t^{-\frac{n}{q}} G_1(x,y,S)+t^{-\frac{n+d}{q}} G_2(x,y,S)+t^{-\frac{p}{q}} G_3(x,y,S), \\
G_1(x,y,S)\equiv \alpha(S) x^2 y, \quad G_2(x,y,S)\equiv \beta(S)y, \quad G_3(x,y,S)\equiv \frac{\gamma(S)y}{\sqrt{x^2+y^2}},\\
\alpha(S)=\alpha_0+\alpha_1 \sin S, \quad 
\beta(S)\equiv\beta_0+\beta_1 \sin S,\quad
\gamma(S)\equiv\gamma_0+\gamma_1 \sin S,
\end{gather*}
$S(t)\equiv \sqrt 2 t+\log t$, $\alpha_j, \beta_j, \gamma_j\in\mathbb R$, and $n,d,p,q\in\mathbb Z_+$ such that $n+d<p$.

We see that system \eqref{Ex2C} describes damped oscillatory perturbations of the pendulum. Let us show that the proposed theory can be applied to this system. 

Note that the corresponding limiting system 
$dx/dt=\partial_y H(x,y)$, $dy/dt=-\partial_x H(x,y)$ with $H(x,y)\equiv y^2/2-\cos x+1$ has a stable equilibrium at the origin. Moreover, the level lines $\{(x,y): H(x,y)=r^2/2\}$ for all $0<|r|<2$, lying in the vicinity of the equilibrium,  correspond to periodic solutions $x_0(t,r)$, $y_0(t,r)$ of the limiting system with a period $T(r)=2\pi/w(r)$, where $w(r)=1-r^2/16+\mathcal O(r^4)$ as $r\to 0$. Define auxiliary $2\pi$-periodic functions 
\begin{gather*}
X(\varphi,r)\equiv x_0\left(\frac{\varphi}{w(r)},r\right), \quad Y(\varphi,r)\equiv y_0\left(\frac{\varphi}{w(r)},r\right).
\end{gather*}
It can easily be checked that 
\begin{gather*}
w(r)\partial_\varphi X(\varphi,r)\equiv Y(\varphi,r), \quad 
w(r)\partial_\varphi Y(\varphi,r) \equiv  -\sin \big(X(\varphi,r)\big), \\
H(X(\varphi,r),Y(\varphi,r))\equiv \frac{r^2}{2},  \quad 
{\hbox{\rm det}}\frac{\partial (X,Y)}{\partial (\varphi,r)}\equiv \begin{vmatrix} \partial_\varphi X & \partial_\varphi Y\\ \partial_r X & \partial_r Y \end{vmatrix} \equiv \frac{r}{w(r)}.
\end{gather*}
Therefore, system \eqref{Ex2C} in the variables $(r,\varphi)$ takes the form \eqref{FulSys} with $\omega(r)\equiv w(r)$, 
\begin{align*}
&f(r,\varphi,S,t)\equiv r^{-1}Y(\varphi,r)G(X(\varphi,r),Y(\varphi,r),S,t), \\ 
&g(r,\varphi,S,t)\equiv - w(r)\partial_r X(\varphi,r)G(X(\varphi,r),Y(\varphi,r),S,t).
\end{align*}
It can easily be checked that these functions satisfy \eqref{FG}:
\begin{align*}
&f(r,\varphi,S,t)\equiv 
t^{-\frac{n}{q}} f_1(r,\varphi,S)+
t^{-\frac{n+d}{q}} f_2(r,\varphi,S)+
t^{-\frac{p}{q}} f_3(r,\varphi,S), \\ 
&g(r,\varphi,S,t)\equiv 
t^{-\frac{n}{q}} g_1(r,\varphi,S)+
t^{-\frac{n+d}{q}} g_2(r,\varphi,S)+
t^{-\frac{p}{q}} g_3(r,\varphi,S),
\end{align*}
where
\begin{align*}
&f_k(r,\varphi,S)\equiv r^{-1}Y(\varphi,r)G_k(X(\varphi,r),Y(\varphi,r),S), \\ 
&g_k(r,\varphi,S)\equiv w(r)\partial_r X(\varphi,r)G_k(X(\varphi,r),Y(\varphi,r),S).
\end{align*}
Since $X(\varphi,r)= r\cos\varphi+\mathcal O(r^3)$, $Y(\varphi,r)= -r\sin\varphi+\mathcal O(r^3)$ as $r\to 0$, we see that assumption \eqref{pFG} holds with $f_{p,0,0}(0)=\gamma_0/2$. Moreover, $S(t)$ has the form \eqref{Sform} with $s_0=\sqrt 2$, $s_1=\dots=s_{q-1}=0$, $s_q=1$, and assumption \eqref{nres} holds.

Let $n=d=1$ and $p=q=3$. Then the change of the variable described in Theorem~\ref{Th1} and constructed in Section~\ref{Sec3} transforms the system into \eqref{FulSys3} with
\begin{align*}
&\Lambda_1(\rho)\equiv \langle f_1(\rho,\varphi,S)\rangle_{\varphi,S}=\frac{\alpha_0 \rho^3}{8}(1+\mathcal O(|\rho|)), \\ 
&\Lambda_2(\rho)\equiv \langle f_2(\rho,\varphi,S)-\mathcal R_2(\rho,\varphi,S)\rangle_{\varphi,S}=\frac{\beta_0 \rho}{2}(1+\mathcal O(|\rho|)), \\ 
&\Lambda_3(\rho)\equiv  \langle f_3(\rho,\varphi,S)-\mathcal R_3(\rho,\varphi,S)\rangle_{\varphi,S}=\frac{\gamma_0}{2}+\mathcal O(|\rho|)
\end{align*} 
as $\rho\to 0$. In this case, $(n+d,q)\in\mathcal Q_3$, $\mu_p=\gamma_0/2$ and assumption \eqref{asnonlin} holds with $m=3$, $\lambda_{n,m}=\alpha_0/8$, $\lambda_{n+d}=\beta_0/2$. Hence, $m_\ast=2$, $\vartheta_m=1/6$ and $\mathcal P(z)\equiv z(\alpha_0 z^2+4\beta_0)/8$. If $\alpha_0<0$ and $\beta_0>0$, then there exists $z_m=\sqrt{-4\beta_0/\alpha_0}$ such that assumption \eqref{Pmas} holds with $\mathcal P(z_m)=0$ and $\mathcal P'(z_m)=\alpha_0 z_m^2/4<0$.
It follows from Theorem~\ref{ThM} that there exists a stable regime for solutions of the system with $r(t)\approx  \varrho_m(t)$  and $\varphi(t)\sim t$, where $\varrho_m(t)$ is a solution of the corresponding truncated equation \eqref{trSys} with asymptotics \eqref{solm}. In particular, $\varrho_m(t)\sim  t^{-1/6}\sqrt{-4\beta_0/\alpha_0}$ as $t\to\infty$. Returning to the variables $(x,y)$, we obtain $H(x(t),y(t))\approx t^{-1/3} (-2\beta_0/\alpha_0)$ (see Fig.~\ref{FigEx2}).
\begin{figure}
\centering
 \includegraphics[width=0.3\linewidth]{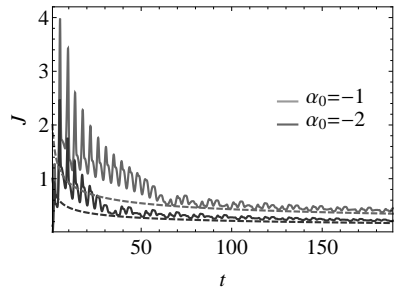}
\caption{\small The evolution of $J(t)\equiv H(x(t),y(t))$  for solutions to system \eqref{Ex2C} with $n=d=1$, $p=q=3$, $\alpha_1=\beta_0=\beta_1=\gamma_0=\gamma_1=1$ and different values of the parameter $\alpha_0$. The dashed curves correspond to $J(t)\equiv 2t^{-1/3}$ and $J(t)\equiv t^{-1/3}$.} \label{FigEx2}
\end{figure}

\section{Conclusion}
Thus, possible long-terms asymptotic regimes of solutions to nonlinear oscillatory systems with damped non-resonant perturbations that do not vanish at the equilibrium have been described. We have shown that, depending on the properties of the perturbations, the asymptotically autonomous system can behave like the corresponding limiting system in the vicinity of the equilibrium. Moreover, new stable regimes associated with the solutions of the corresponding truncated systems may appear. For example, there may be exponentially or polynomially stable solutions with amplitude tending to zero or a non-zero constant. It have been also shown that non-resonant perturbations can lead to instability when the perturbed trajectories leave the neighbourhood of the equilibrium. The results obtained show, in particular, that damped perturbations can be used to control the dynamics of nonlinear oscillatory systems. 

Note that resonant perturbations and higher-dimensional systems have not been considered in this paper. In this case, the proposed theory cannot be applied directly due to the problem of small denominators. These cases deserve special attention and will be discussed elsewhere.

\section*{Acknowledgements}
The research is supported by the Russian Science Foundation (Grant No. 23-11-00009).

}

\begin{thebibliography}{99}

\bibitem{BG08} A. D. Bruno, I. V. Goryuchkina, \textit{Boutroux asymptotic forms of solutions to Painlev\'{e} equations and power geometry}, Doklady Math., 78 (2008), 681--685.

\bibitem {CCT95} C. Castillo-Ch\'{a}vez, H.R. Thieme, \textit{Asymptotically autonomous epidemic models}. In: O. Arino, D. Axelrod, M. Kimmel, M. Langlais (Eds.), Mathematical Population Dynamics: Analysis of Heterogenity, Theory of Epidemics, vol. 1, Wuertz, 1995, p. 33--50.

\bibitem{OS18} O. Sultanov, \textit{Stability and asymptotic analysis of the autoresonant capture in oscillating systems with combined excitation}, SIAM J. Appl. Math., 78 (2018), 3103--3118.

\bibitem{DS22} D. Scarcella, \textit{Weakly asymptotically quasiperiodic solutions for time-dependent Hamiltonians with a view to celestial mechanics}, 2022, arXiv: 2211.06768.

\bibitem{KF13} V. V. Kozlov, S. D. Furta, \textit{Asymptotic Solutions of Strongly Nonlinear Systems of Differential Equations}, Springer, New York, 2013.
 
\bibitem{Pan21} X. Pan, \textit{Stability of smooth solutions for the compressible Euler equations with time-dependent damping and one-side physical vacuum}, Journal of Differential Equations, 278 (2021), 146--188.

\bibitem{Dong22} J. Dong, J. Li, \textit{Analytical solutions to the compressible Euler equations with time-dependent damping and free boundaries}, J. Math. Phys., 63 (2022), 101502.

\bibitem{LRS02} J. A. Langa, J. C. Robinson, A. Su\'{a}rez, \textit{Stability, instability and bifurcation phenomena in nonautonomous differential equations}, Nonlinearity, 15 (2002), 887--903.

\bibitem{KS05} P. E. Kloeden, S. Siegmund, \textit{Bifurcations and continuous transitions of attractors in autonomous and nonautonomous systems}, Internat. J. Bifur. Chaos., 15 (2005), 743--762.

\bibitem{MR08} M. Rasmussen, \textit{Bifurcations of asymptotically autonomous differential equations}, 
Set-Valued Anal., 16 (2008), 821--849.

\bibitem{LM56} L. Markus, \textit{ Asymptotically autonomous differential systems}. In: S. Lefschetz (ed.), Contributions to the theory of nonlinear oscillations III, Ann. Math. Stud., vol. 36, pp. 17--29, Princeton University Press, Princeton, 1956.

\bibitem{LDP74} L. D. Pustyl'nikov, \textit{Stable and oscillating motions in nonautonomous dynamical systems. A generalization of C. L. Siegel's theorem to the nonautonomous case}, Math. USSR-Sbornik, 23 (1974), 382--404.

\bibitem{HT94} H. Thieme, \textit{Asymptotically autonomous differential equations in the plane}, 
Rocky Mountain J. Math., 24 (1994), 351--380.

\bibitem{OS21IJBC} O. A. Sultanov, \textit{Damped perturbations of systems with center-saddle bifurcation}, Internat. J. Bifur. Chaos., 31 (2021), 2150137.

\bibitem{HL75} W. A. Harris and D. A. Lutz, \textit{Asymptotic integration of adiabatic oscillators}, 
J. Math. Anal. Appl., 51 (1975), 76--93.

\bibitem{MP85} M. Pinto, \textit{Asymptotic integration of Second-Order Linear Differential Equations}, J. Math. Anal. Appl., 111 (1985), 388--406.

\bibitem{PN07} P. N. Nesterov , \textit{Averaging method in the asymptotic integration problem for systems with oscillatory-decreasing coefficients}, Differ. Equ. 43 (2007), 745--756.

\bibitem{BN10} V. Burd, P. Nesterov, \textit{Parametric resonance in adiabatic oscillators}, Results. Math., 58 (2010), 1--15.

\bibitem{OS21DCDS} O. A. Sultanov, \textit{Bifurcations in asymptotically autonomous Hamiltonian systems under oscillatory perturbations}, Discrete \& Continuous Dynamical Systems, 41 (2021), 5943--5978.

\bibitem{OS21JMS} O. A. Sultanov, \textit{Decaying oscillatory perturbations of Hamiltonian systems in the plane}, Journal of Mathematical Sciences, 257 (2021), 705--719.

\bibitem{DF78} J. D. Dollard, C. N. Friedman, \textit{Existence of the M{\o}ller wave operators for $V(r)=\frac{\gamma \sin(\mu r^\alpha)}{r^\beta}$}, Annals of Physics, 111 (1978), 251--266.

\bibitem{BD79} M. Ben-Artzi, A. Devinatz, \textit{Spectral and scattering theory for the adiabatic oscillator and related potentials}, J. Math. Phys., 20 (1979), 594--607. 

\bibitem{OS22JMS} O. A. Sultanov, \textit{Capture into resonance in nonlinear oscillatory systems with decaying perturbations}, Journal of Mathematical Sciences, 262 (2022), 374--389.

\bibitem{OS23DCDSB} O. A. Sultanov, \textit{Resonances in asymptotically autonomous systems with a decaying chirped-frequency excitation}, Discrete and Continuous Dynamical Systems --- B, 28 (2023), 1719--1749.

\bibitem{WW66} W. Wasow, \textit{Asymptotic Expansions for Ordinary Differential Equations}, John Wiley and Sons, Inc., New York, 1966. 

\bibitem{BM61} N. N. Bogolubov, Yu. A. Mitropolsky, \textit{Asymptotic Methods in Theory of Non-linear Oscillations}, Gordon and Breach, New York, 1961.

\bibitem{AN84} A. I. Neishtadt, \textit{The separation of motions in systems with rapidly rotating phase}, J. Appl. Math. Mech., 48 (1984), 133--139.

\bibitem{AKN06} V. I. Arnold, V. V. Kozlov, A. I. Neishtadt, \textit{Mathematical Aspects of Classical and Celestial Mechanics}, Springer, Berlin, 2006.

\bibitem{LK15} L. A. Kalyakin, \textit{Lyapunov functions in theorems of justification of asymptotics}, Mat. Notes,
98 (2015), 752--764.

\end{thebibliography}
\end{document}